\documentclass[12pt,a4paper]{article}

\usepackage[brazil]{babel}
\usepackage[hmargin=3cm,vmargin=3cm,bmargin=3cm]{geometry}

\usepackage{amssymb,amsmath,amsthm,amsfonts,graphicx}
\usepackage{enumerate}
\usepackage[utf8]{inputenc}
\usepackage{array} 
\usepackage{multicol}
\usepackage{float} 
 \usepackage{booktabs}
  \usepackage[dvipsnames]{xcolor}
 \usepackage{xcolor}
 \usepackage{tcolorbox}
 \usepackage{listings}
 \usepackage{multicol}
\usepackage{array, ltablex}
\usepackage{enumitem}%
\usepackage{makecell, etoolbox}
\setcellgapes{\itemsep}
\AtBeginEnvironment{tabularx}{\setcellgapes[b]{\itemsep}\makegapedcells}
\AtEndEnvironment{tabularx}{\vskip\dimexpr-\itemsep-\topsep\relax}
\newcolumntype{I}{@{\hskip\leftmargin} >{\llap{\textbullet\hskip\labelsep}}X}

\newtheorem{theorem}{Teorema}[section]
\newtheorem{proposicao}{Proposição}

\newtheorem{exemplo}{Exemplo}[section]
\newtheorem{definition}{Definição}[section]

\newtheorem{observacao}{Observação}

\begin{document}
\def\n{\noindent}
\begin{center}
\textbf{Formas Quadráticas e  Método de Fatorização de Gauss.}

\vspace{0.5 cm}
Seminário em Matemática e Aplicações.\\

\vspace{0.2 cm}
Mónica Celis.\\
\vspace{0.2 cm}
Paulo Almeida.\\

\end{center}
\vspace*{2pt}
\section{ Alguns resultados sobre formas quadráticas. }

\subsection{Formas quadráticas.}

\begin{definition}
Uma \textbf{forma quadrática} é uma função em duas indeterminadas da forma:
\[ax^2+2bxy+cy^2\]
onde a,b,c são números inteiros dados. \\

\n Denotaremos a forma $ax^2+2bxy+cy^2$ por $(a,b,c).$
\end{definition}

\begin{definition}

O \textbf{determinante} da forma $(a,b,c)$ será definido como o número $b^2-ac.$
\end{definition}

\begin{theorem}
Se um número $M$ pode-se representar  pela forma $(a,b,c)$  tal que os valores das indeterminadas que fazem que isto aconteça, são relativamente primos entre si, $b^2-ac$ será um resíduo quadrático do número $M$.

\end{theorem}
\n \textbf{Demonstração.}\\

\n Sejam $m,n$ os valores das indeterminadas tais que

\[M= am^2+2bmn+cn^2\]

\n e sejam $\mu\ e\  \nu$ tais que $\mu m+\nu n =1$ pois $m,n$ são relativamente primos, então\\

\n $(am^2+2bmn+cn^2)(a\nu^2-2b\mu \nu+c \mu^2)  $\\
$=a^2m^2 \nu^2-2abm^2 \mu \nu+acm^2\mu^2+2abmn\nu^2-4b^2mn\mu\nu+2bcmn\mu^2+acn^2\nu^2-2bcn^2\mu\nu+c^2n^2\mu^2$\\
\n $=$\textcolor{blue} {$b^2m^2\mu^2$} $+ 2bcmn\mu^2 +c^2n^2\mu^2  -2abm^2 \mu \nu-2b^2mn\mu\nu$ \textcolor{red}{$-2acmn\mu\nu$} $ -2bcn^2\mu \nu + a^2m^2 \nu^2 +2abmn\nu^2 $  \textcolor{olive}{$+b^2n^2\nu^2$}  \textcolor{blue} {$-b^2m^2\mu^2$}  $-2b^2mn\mu\nu $ \textcolor{olive}{$-b^2n^2\nu^2$} $+acm^2\mu^2 $  \textcolor{red}{$+2acmn\mu\nu$}  $+acn^2\nu^2$ \\
$=( bm\mu+ cn\mu)^2 - 2 \mu \nu(bm+cn)(am+bn) +(am\nu +bn\nu)^2 -b^2(m^2\mu^2+2mn\mu\nu +n^2\nu^2) +ac(m^2\mu^2+2mn\mu\nu +n^2\nu^2)$\\
$=\left[ \mu(bm+cn)- \nu (am+bn)\right]^2 -(b^2-ac)(m\mu+n\nu)^2. $

\n Assim,
\[M(a\nu^2-2b\mu \nu+c \mu^2)= \left[ \mu(bm+cn)- \nu (am+bn)\right]^2 -(b^2-ac)(1)^2,\]

\n portanto, 
\[(b^2-ac) \equiv   \left[ \mu(bm+cn)- \nu (am+bn)\right]^2 \ (mod\ M).\]

\n Logo, $b^2-ac$ é resíduo quadrático mod $M$. \ \ \ \ \ \ \ \ \ \ \ \  \ \ \ \ \ \ \ \ \ \ \ \ \ \  \ \ \ \ \ \ \ \ \ \ \ \ \ \  \ \ \ \ \ \ \ \ \ \ \ \ \  \ \ \ $\square$\\

\n Do teorema anterior temos que 

\[  \mu(bm+cn)- \nu (am+bn) \]

\n é o valor da expressão $\sqrt{(b^2-ac)} \ mod\ M. $ 
\begin{observacao}
O valor de $\sqrt{(b^2-ac)} \ mod\ M$ não muda se os valores de $\mu$ e $\nu$ mudam, isto é, se $\mu m+\nu n =1$ e se também $\mu' m+\nu' n =1$ então,  considerando
\[ \mu(bm+cn)- \nu (am+bn) = \vartheta,\ \ \ \   \mu'(bm+cn)- \nu' (am+bn) = \vartheta',\]


\n se multiplicamos a equação $\mu m+\nu n =1$ por $\mu'$ e $\mu' m+\nu' n =1$ por $\mu$, e subtraímos obtemos  $\mu'-\mu= n(\mu' \nu-\mu\nu'). $ Se multiplicamos a equação $\mu m+\nu n =1$ por $\nu'$ e $\mu' m+\nu' n =1$ por $\nu$, e subtraímos obtemos  $\nu'-\nu= m(\mu \nu'-\mu' \nu). $ Daqui resulta que, 

\[\begin{array}{ll}

\vartheta'-\vartheta&= (\mu'-\mu)(bm+cn)-(\nu'-\nu)(am+bn)\\
&= n(\mu' \nu-\mu\nu')(bm+cn)- m(\mu \nu'-\mu' \nu)(am+bn)\\
&= n(\mu' \nu-\mu\nu')(bm+cn)+ m(\mu' \nu-\mu \nu')(am+bn)\\
&= (\mu' \nu-\mu\nu')\left[ n(bm+cn)+ m(am+bn)\right] \\
&= (\mu' \nu-\mu\nu')\left[ am^2+2bmn+cn^2\right] \\
&= (\mu' \nu-\mu\nu')M, \\
\end{array}\]

\n  ou seja, $\vartheta' \equiv \vartheta$ (mod $M).$

\n Portanto, de qualquer forma que sejam determinados  o $\mu$ e o $\nu$, o valor de $\mu(bm+cn)- \nu (am+bn)$ não pode apresentar valores diferentes (i.e. não congruentes)  da expressão $\sqrt{(b^2-ac)} \ mod\ M. $ \\

 \end{observacao}
\n Se $\vartheta$ é um valor  de $\mu(bm+cn)- \nu (am+bn)$, dizemos que a representação do número  $M$ pela forma $ax^2+2bxy+cy^2$  onde $x=m\ e\ y=n$ \textbf{pertence ao valor $\vartheta$} de $\sqrt{(b^2-ac)} \ mod\ M.$\\

\begin{definition}
Se uma forma $F$ com  indeterminadas $x,y$  pode ser transformada em outra forma $F'$ com indeterminadas $x',y'$ mediante as substituições

\[x=\alpha x' +\beta y', \ \ \ \ \ \ \ \ \  y=\gamma x'+\delta y' \]

\n onde $\alpha, \beta, \gamma \ e\ \delta$ são inteiros, dizemos que a forma $F$ \textbf{implica} $F'$ ou que $F'$ está \textbf{contida} em $F$.
\end{definition}

\n Sejam $F$ a forma $ ax^2+2bxy+cy^2$ e $F'$ a forma  $a'x'^2+2b'x'y'+c'y'^2$ da definição anterior. Então 
\[a'x'^2+2b'x'y'+c'y'^2\ = a(\alpha x' +\beta y')^2+2b(\alpha x' +\beta y')(\gamma x'+\delta y')+c(\gamma x'+\delta y')^2. \]
$a'x'^2+2b'x'y'+c'y'^2 = (a \alpha^2+2b\alpha \gamma + c \gamma^2) x'^2  +2(a \alpha \beta +b(\alpha \beta+\beta \gamma) + c \gamma \delta )x'y' +$
\hspace*{4cm}$+(a \beta^2+2b\beta \delta +c \delta^2) y'^2.$

\n Assim, obtemos as seguintes equações:\\

\n $a'= a \alpha^2+2b\alpha \gamma + c \gamma^2$.\\
$b'= a \alpha \beta +b(\alpha \beta+\beta \gamma)+ c \gamma \delta $.\\
$c'= a \beta^2+2b\beta \delta +c \delta^2.$\\

\n Logo, 
\begin{equation} \label{eq1}
b'^2-a'c'=(b^2-ac)(\alpha\delta-\beta\gamma)^2
\end{equation}

\n De  \eqref{eq1} podemos concluir que o determinante de $F'$  é divisível pelo determinante de $F$ e o quociente entre eles é um quadrado, logo os determinantes terão o mesmo sinal. 

\begin{observacao}
Se $F'$ está contida em $F$ e $F$ está contida em $F'$, então os determinantes das formas serão iguais e $(\alpha\delta -\beta\gamma)^2=1$. Neste caso, diremos que $F$ e $F'$ são \textbf{equivalentes}.
\end{observacao}

\begin{definition}
A substituição $x=\alpha x' +\beta y',  y=\gamma x'+\delta y'$ será chamada de transformação\textbf{ própria} se  $(\alpha\delta -\beta\gamma)$ é um número positivo e \textbf{imprópria} se $(\alpha\delta -\beta\gamma)$ é um número negativo. 
\end{definition}
\begin{definition}
A forma $F'$ está \textbf{contida propriamente ou imprópria\-mente} na forma $F$, se $F$ pode ser transformada  na forma $F'$ por uma transformação própria ou imprópria, respetivamente.
\n Se $F'$ está contida propriamente  em $F$ e  se $F$ está contida propriamente em $F'$  então dizemos que as formas são \textbf{propriamente equivalentes}.  Se $F'$ está contida impropriamente  em $F$ e  se $F$ está contida impropriamente em $F'$  então dizemos que as formas são \textbf{impropriamente equivalentes}.
\end{definition}
\begin{theorem} \label{T2}
Se a forma $F$ implica a forma $F'$ e a forma $F'$ implica a forma $F''$, então a forma $F$ também implica a forma $F''$.
\end{theorem}
\begin{observacao}
Sabemos que se as formas $F,F'$ são equivalentes então  $(\alpha\delta -\beta\gamma)^2=1$, se a transformação é própria então $\alpha\delta -\beta\gamma=1$; se a transformação é imprópria então $\alpha\delta -\beta\gamma=-1.$ 
\end{observacao}
\begin{exemplo}
A forma $-13x_1^2-12x_1y_1-2y_!^2$ está contida propriamente na forma $2x^2-8xy+3y^2$, a transformação própria usada é $x=2x_1+y_1, y= 3x_1+2y_1$  $(\alpha\delta -\beta\gamma= 2\cdot 2 - 1 \cdot 3 = 1)$ e a forma $2x^2-8xy+3y^2$ está contida propriamente na forma  
 $-13x_1^2-12x_1y_1-2y_1^2$, a transformação usada é  $x=2x_1-y_1, y= -3x_1+2y_1$. Assim, as formas $2x^2-8xy+3y^2$ e $-13x'^2-12x_1y_1-2y'^2$ são propriamente equivalentes. 
\end{exemplo}
\n A forma $(c,-b,a)$ é propriamente equivalente a forma $(a,b,c)$  e as formas $(a,-b,c)$ e $(c,b,a)$ são impropriamente equivalentes a forma $(a,b,c)$.\\

\n $ax^2+2bxy+cy^2 $ $\xrightarrow[x=0x_1-y_1, y=x_1+0y_1 ]{} cx_1^2-2bx_1y_1+ay_1^2.$ \\

\n   $ax^2+2bxy+cy^2  $ $\xrightarrow[x=x_1+0y_1, y=0x_1-y_1 ]{} ax_1^2-2bx_1y_1+cy_1^2.$\\

\n $ax^2+2bxy+cy^2  $ $\xrightarrow[x=0x_1+y_1, y=x_1+0y_1]{} cx_1^2+2bx_1y_1+ay_1^2.$\\

\n Pelo anterior, qualquer forma equivalente à forma $(a,b,c)$ é propriamente equivalente a ela mesma ou à forma $(a,-b,c)$. Similarmente, se qualquer forma implica a forma $(a,b,c)$  ou está contida nesta, a forma será propria\-mente equivalente à forma $(a,b,c)$  ou à forma $(a,-b,c)$ ou estará contida propriamente numa delas. Chamamos às formas $(a,b,c)$ e  $(a,-b,c)$ \textbf{formas opostas}.


\begin{definition}
Se as formas $(a,b,c)$ e $(a',b',c')$  possuem o mesmo determinante e além disso $c=a'$ e $b \equiv -b'\ (mod\ c)$ então dizemos que as formas são \textbf{contíguas} ou \textbf{vizinhas.}  
\end{definition}

\begin{observacao} Formas contíguas são sempre propriamente equivalentes. 
\end{observacao}

\begin{theorem}
Se a forma $F$ implica a forma $F',$ qualquer número que se possa representar por $F'$ também pode ser representado por $F$.
\end{theorem}


\n \textbf{Demonstação.}

\n Sejam $x\ e\ y$, $x'\ e\ y'$ as indeterminadas das formas $F$ e $F'$ respetivamente, e suponhamos que o número $M$ se represente por $F'$ fazendo  $x'=m$ e $y'=n$.\\
\n A forma $F$ se transforma na forma $F'$ pela substituição:
\[x=\alpha x' +\beta y', \ \ \ \ \ \ \ \ \  y=\gamma x'+\delta y' \]
\n Então, o número $M$ pode ser representado pela forma $F$ se consideramos 
\[x=\alpha m + \beta n, \ \ \ \ \ \ \ \ \  y=\gamma m + \delta n.  \] \begin{flushright}
$\square $
\end{flushright}
\begin{observacao}
Se $F$ e $F'$ são equivalentes e se o número $M$ pode-se representar por uma das formas então também pode-se representar pela outra.
\end{observacao}
\begin{theorem}
Se as formas 
\[F:ax^2 +2bxy+cy^2\]
\[F':ax'^2 +2bxy+cy'^2\]
\n são equivalentes, se o determinante é igual a $D$, e se a forma $F'$ é transformada na forma $F$  pela substituição
\[ x'= \alpha x + \beta y, \ \ \ \ \  y'= \gamma x + \delta y. \]

\n Se além disso, o número $M$ pode-se representar por $F$ sendo $x=m$ e $y=n$ e por $F'$ sendo 

\[x'= \alpha m+ \beta n= m', \ \ \ \ \ y'= \gamma m + \eta n = n'\]

\n de maneira que $m$ seja relativamente primo a $n$ e portanto também $m'$ e $n'$ sejam relativamente primos, então as duas representações pertencerão ou ao mesmo valor da expressão $\sqrt{D}\ mod\ M$ ou a valores opostos de acordo como a transformação da forma $F'$ em $F$ seja própria ou imprópria.    
\end{theorem}

\subsection{Formas quadráticas com determinante negativo.}

\begin{theorem}
Se o número $M$ pode-se representar pela forma $ax^2+2bxy+cy^2$, designando os valores $x=m$ e $y=n$ onde $(m,n)=1$ e se o valor da expressão $\sqrt{D}\ mod\ M$  ao qual pertence está representação é $N$, então as formas $(a,b,c)$ e $(M,N,(N^2-D)/M)$ serão propriamente equivalentes. 
\end{theorem}
\begin{definition}
Uma forma $(A,B,C)$ cujo determinante é igual a $-D$ onde $D$ é um número positivo e na qual $A$ é menor que $\sqrt{4D/3}$ ou $C$ e  é maior que $2B$ é chamada de \textbf{forma reduzida.}
\end{definition}
\begin{theorem}
Dada uma forma qualquer $(a,b,a_1)$ cujo determinante negativo é igual a $-D$ onde $D$ é um número positivo,  então pode-se encontrar uma forma reduzida $(A,B,C)$ propriamente equivalente a $(a,b,a_1).$
\end{theorem}
\n \textbf{Demonstração.}

\n Sejam $b_1$ o resíduo absolutamente mínimo\footnote{Todo número terá um resíduo mod $m$, tanto na sucessão $0,1,2,\cdots,m-1$,
como na sucessão $ 0,-1,-2,\cdots,-(m-1)$ chamamos a estes resíduos de \textbf{\textit{resíduos mínimos}.} A menos que 0 seja um resíduo teremos sempre um resíduo positivo e um negativo, destes dois resíduos o que seja menor que $\frac{m}{2}$ em modulo chamaremos de \textbf{\textit{resíduo absolutamente mínimo. }}Por exemplo, $-12 \equiv 3\ mod\ 5 $, 3 é  resíduo mínimo positivo, $-12 \equiv -2\ mod\ 5 $, -2 é  resíduo mínimo negativo. Logo, $-2$ é o resíduo absolutamente mínimo.   } do número $-b$ módulo  $a_1$ e $a_2= \dfrac{b_1^2+D}{a_1}. $

\n $a_2$ é inteiro, pois $b_1^2 \equiv b^2$ e $b_1^2 +D \equiv b^2+D = b^2-b^2+aa_1= aa_1 \equiv 0\ (mod\ a_1).$\\ 
\n Se $a_2 < a_1$, seja $b_2$  o  resíduo absolutamente mínimo do número $-b_1$ módulo $a_2$ e  $a_3= \dfrac{b_2^2+D}{a_2}. $\\

\n Se $a_3< a_2$, seja $b_3$  o resíduo absolutamente mínimo do número $-b_2$ módulo $a_3$ e  $a_4= \dfrac{b_3^2+D}{a_3}. $

\n Continuar está operação até chegar na progressão $a_1,a_2,a_3,\cdots$ a um termo $a_{m+1}$ o qual não será menor que seu antecessor $a_m$. Isto deve acontecer pois caso contrário ter-se-ia uma progressão infinita decrescente de números inteiros. Assim a forma $(a_m,b_m,a_{m+1})$ satisfará as condições.
 
\n Na progressão das formas  $(a,b,a_1),(a_1,b_1,a_2),(a_2,b_2,a_3), \cdots $ cada uma é contígua à sua antecessora, pelo qual a ultima será equivalente à primeira (teorema \ref{T2}). \\
\n Como $b_m$ é o menor resíduo absoluto de $-b_{m-1}$ módulo $a_m$ não será maior que $\frac{1}{2} a_m.$\\
\n Como $a_ma_{m+1}= D+b_{m}^2$ e $a_{m+1}$ é maior que $a_m$ então  $a_m^2 \leq D+b_{m}^2$  e como $b_m$  é menor que $\frac{1}{2} a_m$ então  $a_m^2$ também será menor que $D+ \frac{1}{4} a_m^2$. Logo, $\frac{3}{4} a_m^2 $  é menor ou igual  que $D$. Assim, $a_m$  é menor que $\sqrt{\frac{4}{3}D.}$ \ \ \ \ \ \ \ \ \  \ \ \ \ \ \ \ \ \  \ \ \ \ \ \ \ \ \  \ \ \ \ \ \ \ \ \  \ \ \ \ \ \ \ \ \  \ \ \ \ \ \ \ \ \  \ \ \ \ \ \ \ \ \  \  $ \square$

\begin{exemplo}
 Determinar uma forma reduzida propriamente equivalente com  a forma $(304,217,155).$\\
 
 O determinante da forma é $-31$, $D=31$.
 
 \begin{itemize}
 \item[•] $(a,b,a_1)= (304,217,155)$
 
 $-217 \equiv -62\ mod\ 155 \Rightarrow b_1=-62$
 
 $a_2= \frac{b_1^2+D}{a'}= \frac{(-62)^2+31}{155}=\frac{3875}{155}=25 \Rightarrow a_2=25$ 
 
 $a_2< a_1$
 \item[•] $(a_1,b_1,a_2)= (155,-62,25)$
 
  $62\equiv 12\ mod\ 25 \Rightarrow b_2=12$
  
   $a_3= \frac{b_2^2+D}{a_2}= \frac{(12)^2+31}{25}=7 \Rightarrow a_3=7$ 
 
 $a_3< a_2$
 
  \item[•] $(a_2,b_2,a_3)= (25,12,7)$
 
  $-12\equiv 2\ mod\ 7\Rightarrow b_3=2$
  
   $a_{4}= \frac{b_3^2+D}{a_3}= \frac{(2)^2+31}{7}= 5\Rightarrow a_{4}=5$ 
 
 $a_{4}< a_3$
 
  \item[•] $(a_3,b_3,a_{4})= (7,2,5)$
 
  $-2\equiv -2\ mod\ 5\Rightarrow b_{4}=-2$
  
   $a_{5}= \frac{(b_{4})^2+D}{a_{4}}= \frac{(-2)^2+31}{5}= 7\Rightarrow a_{5}=7$ 
 
\n Como $a_{5} \geq a_{4}$ então a forma buscada é $(5,-2,7).$
 \end{itemize}

\end{exemplo}
\n \textbf{Métodos para encontrar todas as formas reduzidas de determinante $-D$} \footnote{o número de formas reduzidas que tem como determinante $-D$ é finito e relativamente pequeno em relação ao número $D$. } \\ 

\n Designamos por $(a,b,c) $ as formas reduzidas com determinante $-D$ que queremos encontrar. \\

\n  \underline{\textit{Primeiro Método:}} Tome-se  para $a$ todos os valores positivos e negativos não maiores que $\sqrt{\frac{4}{3}D}$ para os quais $-D$ é resíduo quadrático. Para cada um desses $a$, seja $b$ o valor da expressão $\sqrt{-D}\ mod\ a$ que não seja maior que $\frac{a}{2}$, tanto valores positivos como negativos podem ser usados. \\
Para cada par de valores de $a$ e $b $, fazemos $c=\frac{D+b^2}{a}$. \\
Se deste método resultam formas onde $c<a$ então ditas formas devem ser descartadas. \\

\n  \underline{\textit{Segundo Método:}}  Tome-se  para $b$ todos os números positivos e negativos  os quais em valor absoluto são menores  que  $\sqrt{\frac{D}{3}}$. Para cada $b$, decompor
$b^2+D$  de todas as maneiras em dois pares de fatores onde nenhum dos dois sejam  menores que $2b$ (valores positivos e negativos devem ser usados). Quando os fatores são diferentes, consideramos como $a$ o menor deles e o outro como $c$. Da forma como foi construido $a$ sabemos que $a$ não é maior que $\sqrt{\frac{4}{3}D}.$
\begin{exemplo}
Encontrar todas as formas reduzidas com determinante igual a $-85.$\\

\n Seja $D=85$.\\

\n \textit{Primeiro Método}\\

\n O limite para os valores de $a$ é $\sqrt{\frac{340}{3}} \thickapprox 10,64$

\n Assim, os números entre 0 e 10 que são resíduos quadráticos módulo $-85$ são: 1,2,5,10. As formas são as seguintes:

\begin{itemize}
\item[•] $(1,0,85)$, $a=1, b=\sqrt{-85}\ mod\  1 \Rightarrow b=0, c= \frac{85}{1}=85.$

\item[•] $(2,1,43)$, $a=1, b=\sqrt{-85}\ mod\  2 \Rightarrow b=1, c= \frac{85+1}{2}=43.$

\item[•]  $(2,-1,43)$, $a=1, b=\sqrt{-85}\ mod\  2 \Rightarrow b=-1, c= \frac{85+1}{2}=43.$

 \item[•]  $(5,0,17)$, $a=5, b=\sqrt{-85}\ mod\  5 \Rightarrow b=0, c= \frac{85}{5}=17.$
 
 \item[•]  $(10,5,11)$, $a=10, b=\sqrt{-85}\ mod\  10 \Rightarrow b=5, c= \frac{85+25}{10}=11.$
 
  \item[•]  $(10,-5,11)$, $a=10, b=\sqrt{-85}\ mod\  10 \Rightarrow b=-5, c= \frac{85+25}{10}=11.$
  
  \n Considerando os valores negativos temos  as formas:
 
  \item[•]  $(-1,0,-85), (-2,1,-43),(-2,-1-43),(-5,0,-17),(-10,5,-11),$\\
  $(-10,-5,-11).$
\end{itemize}

\n Logo, as formas reduzidas de determinante $-85 $ são:
\[(1,0,85),(2,1,43),(2,-1,43), (5,0,17),(10,5,11),(10,-5,11),(-2,1,-43),\]
\[ (-1,0,-85),(-2,-1-43),(-5,0,-17),(-10,5,-11),(-10,-5,-11).\]

\vspace{.2cm}

\n \textit{Segundo Método}\\

O limite dos valores de $b$ é $\sqrt{\frac{D}{3}}=\sqrt{\frac{85}{3}} \thickapprox 5,32.$

\begin{itemize}
\item[•] $b=0$\\
$b^2+D=85=1\cdot 85 \longrightarrow a=1, c=85 \longrightarrow(1,0,85), (-1,0,-85).$ \\
$b^2+D=85=5\cdot 17 \longrightarrow a=5, c=17 \longrightarrow(5,0,17),(-5,0,-17).$ 
\item[•] $b=1$\\
$b^2+D=1+85=86=2\cdot 43 \longrightarrow a=2, c=43 \longrightarrow(2,1,43),(-2,1,-43).$ 
\item[•] $b=-1 \longrightarrow (2,-1,43),(-2,-1,-43).$

\item[•] $b=2, (b=-2).$\\
$b^2+D=4+85=89$, $89$ não pode-se escrever como produto de dois fatores, onde cada fator não seja menor que $4$ ($=2b$). Logo, não existe formas para este caso. 

\item[•] $b=3,  (b=-3).$\\
$b^2+D=9+85=94=2\cdot 47$, $94$ não pode-se escrever como produto de dois fatores, onde cada fator não seja menor que $6$ ($=2b$). 

\item[•] $b=4,  (b=-4).$\\
$b^2+D=16+85=101$, $101$ não pode-se escrever como produto de dois fatores, onde cada fator não seja menor que $8$ ($=2b$). 

\item[•] $b=5.$\\
$b^2+D=25+85=110=10\cdot 11 \longrightarrow a=10, c=11 \longrightarrow$\\
\hspace*{7.3cm} $\longrightarrow(10,5,11),(-10,5,-11).$ 

\item[•] $b=-5 \longrightarrow (10,-5,11),(-10,-5,-11).$

\end{itemize}

\n Assim, as formas reduzidas  com determinante igual a $-85$ são:
\[(1,0,85),(2,1,43),(2,-1,43), (5,0,17),(10,5,11),(10,-5,11),\]
\[(-2,1,-43), (-1,0,-85),(-2,-1-43),(-5,0,-17),(-10,5,-11),\]
\[(-10,-5,-11).\]
\end{exemplo}

\begin{observacao}
Se entre as formas reduzidas de um determinante dado, descartamos uma ou outra das formas não iguais que são propriamente equivalentes, as formas restantes terão a seguinte propriedade: qualquer forma com o mesmo determinante será propriamente equivalente a uma e somente uma delas. De onde, todas as formas do mesmo determinante podem ser distribuídas  em tantas classes como formas tenham ficado. Por exemplo, para $D=85$ as formas que ficam são:

\begin{equation}\label{eq2}
\begin{split}
(1,0,85),(2,1,43), (5,0,17),(10,5,11),\\
(-1,0,-85),(-2,1-43),(-5,0,-17),(-10,5,-11)
\end{split}
\end{equation}

\n pois $(2,-1,43)$ é propriamente equivalente a $(2,1,43)$, $(10,-5,11)$ é propriamente equivalente a $(10,5,11)$,$ (-2,-1-43)$  é propriamente equivalente a $(-2,1-43)$ e $(-10,-5,-11)$ é propriamente equivalente a $(-10,5,-11).$\\

\n Assim, as formas de determinante $-85$ podem ser distribuídas em oito classes, conforme elas sejam propriamente equivalentes a alguma das formas de \eqref{eq2}.
\end{observacao}

\subsection{Formas quadráticas com determinante positivo não quadrado.}

\begin{definition}
Uma forma $(A,B,C)$ com determinante positivo não quadrado $D$ tal que o valor de $A$ tomado positivamente  está entre $\sqrt{D}+B$ e $\sqrt{D}-B$, onde $B$ é positivo e menor que $\sqrt{D}$ será chamada de \textbf{forma reduzida.}
\end{definition}


\begin{theorem}
Dada uma forma qualquer $(a,b,a_1)$ cujo determinante positivo não quadrado é igual a $D$,  então pode-se encontrar uma forma reduzida $(A,B,C)$ propriamente equivalente a $(a,b,a_1).$
\end{theorem}
\n \textbf{Demonstração.}\\

\n Seja $b_1\equiv -b\ mod\ a_1$ de modo que esteja entre  $\sqrt{D}$ e $\sqrt{D}\mp a_1$ (considerando o sinal superior quando  $a_1$ seja positivo e o inferior quando $a'$ seja negativo).\\
\n Seja $a_2= \frac{b_1^2 -D}{a_1}$, o qual é inteiro pois $b_1^2 -D\equiv b^2-D\equiv aa_1\equiv 0\ mod\ a_1.$\\
\n Agora, se $a_2<a_1$ então consideraremos $b_2 \equiv -b_1\ mod\ a_2$ situado entre $\sqrt{D}$ e $\sqrt{D}\mp a_2$ (segundo $a_2$ seja positivo ou negativo) e $a_3= \frac{b_2^2 -D}{a_2}$.\\
\n Se novamente $a_3<a_2$, seja $b_3 \equiv -b_2\ mod\ a_3$ situado entre $\sqrt{D}$ e $\sqrt{D}\mp a_3$  e $a_{4}= \frac{b_3^2 -D}{a_3}$. O processo deve continuar até chegar na progressão $a_1,a_2,a_2,\cdots$ a um termo $a_{m+1}$ o qual não será menor que $a_m$.\\
\n Verifiquemos que a  forma $(A,B,C)=(a_m,b_m,a_{m+1})$ satisfaz todas as condições.

\begin{itemize} 

\item[•]  Na progressão das formas  $(a,b,a_1),(a_1,b_1,a_2),(a_2,b_2,a_3), \cdots $ cada uma é contí\-gua a sua antecessora, pelo qual a última será equivalente à primeira (Teorema \ref{T2}), isto é, $(A,B,C)=(a_m,b_m,a_{m+1})$ será equivalente a $(a,b,a_1)$.
\item[•]$B<\sqrt{D} \ e\  positivo.$\\
 Como $B$ está situado entre  $\sqrt{D}$ e $\sqrt{D}\mp A$ (considerando o sinal superior quando A é positivo e o sinal inferior quando A é negativo), se consideramos $p=\sqrt{D}-B,$ $q= B- (\sqrt{D}\mp A)$, $p,q$ são positivos. Agora, pode-se verificar que $q^2+2pq+2p\sqrt{D}= D+A^2-B^2$, daqui que $D+A^2-B^2$ será um número positivo, seja $r=D+A^2-B^2$\\
Como $D=B^2-AC$ então $r= A^2-AC$ assim $A^2-AC$ será um número positivo mas dado que $A$ é menor que $C$ então devemos ter que $AC$ é negativo pois caso contrario  $A^2-AC$ não seria positivo. Logo $A $ e $ C$ devem ter sinais contrários.  Portanto, $B^2=D+AC<D$ o que implica que $B<\sqrt{D}.$\\ 
\n  Além disso,  como $-AC=D-B^2$, então $AC<D$ e daqui $AA<AC<D$, logo $A<\sqrt{D}$. Pelo que, $\sqrt{D} \mp A$ será positivo e consequentemente também $B$ será positivo, pois situa-se entre $\sqrt{D}$ e $\sqrt{D}\mp A$.

\item[•] A está situado entre $\sqrt{D}+B$ e $\sqrt{D}-B$.

Como $B$ e $\sqrt{D} \mp A$ são positivos então $\sqrt{D} \mp A +B$ é positivo e como $-q= \sqrt{D} \pm A -B$ é negativo (pois $q$ é positivo),  então $ \pm A $ deverá estar situado entre $\sqrt{D}+B$ e $\sqrt{D}-B$.

\end{itemize}
\begin{exemplo}
Encontrar uma forma reduzida  propriamente equivalente a $(67,97,140).$\\
\n O determinante da forma é  $97^2-67\cdot 140=29.$

\begin{itemize}
\item[•]  $(a,b,a_1)=(67,97,140)$

$b_1 \equiv -97\ mod\ 140 \Rightarrow b_1=-97$\\
 $(b_1$ deve estar situado entre $\sqrt{D} -140 \approx -134.61 $ e  $\sqrt{D} \approx 5,38).$  \\
 $a_2= \frac{(-97)^2-29}{140}=67.$\\
 $a_2<a_1$
 
 \item[•]  $(a_1,b_1,a_2)=(140, -97,67)$
 
 $b_2 \equiv 97\ mod\ 67 \Rightarrow b_2=-37$\\
 $(b_2$ deve estar situado entre $\sqrt{D} -67 \approx -61,61$ e  $\sqrt{D} \approx 5,38).$  \\
 $a_3= \frac{(-37)^2-29}{67}=20.$\\
 $a_3<a_2$
 
  \item[•]  $(a_2,b_2,a_3)=(67,-37,20)$
 
 $b_3 \equiv 37\ mod\ 20 \Rightarrow b_3=-3$\\
 $(b_3$ deve estar situado entre $\sqrt{D} -20 \approx -14,61$ e  $\sqrt{D} \approx 5,38).$  \\
 $a_{4}= \frac{(-3)^2-29}{20}=-1.$\\
 $a_{4}<a_3.$
  
   \item[•]  $(a_3,b_3,a_{4})=(20,-3,-1).$
 
 $b_{4} \equiv 3\ mod\ 1\Rightarrow b_{4}=5$\\
 $(b_4$ deve estar situado entre  $\sqrt{D} -1 \approx 4,38$ e $\sqrt{D} \approx 5,38)$.  \\
 $a_{5}= \frac{(5)^2-29}{-1}=4.$\\
 $a_{5}\nless a_{4}.$
 
 \n Logo, a forma buscada é $(-1,5,4).$
 
\end{itemize}

\end{exemplo}
\n \textbf{Métodos para encontrar todas as formas reduzidas de determinante positivo não quadrado $D$.}\\

\n Designamos por $(a,b,c) $ as formas reduzidas com determinante $D$ que queremos encontrar. \\

\n  \underline{\textit{Primeiro Método:}} Tome-se para $a$  todo os números (tanto positivos como negativos) que em módulo são menores que $2\sqrt{D}$ dos quais $D$ é resíduo quadrático.\\
\n Para cada $a$ em particular, seja $b$ igual a todos os valores positivos da expressão $\sqrt{D}\ mod\ a$ situados entre $\sqrt{D}$ e $\sqrt{D}\mp a$(considerando o signo superior quando  $a$ seja positivo e o inferior quando $a$ seja negativo). Para cada um dos valores de $a$ e $b$ seja $c=\frac{b^2-D}{a}$\\
\n Deve-se descartar as formas onde $\pm a$ não esteja entre $\sqrt{D}+b$ e $\sqrt{D}-b.$\\

\begin{exemplo}
Encontrar todas as formas reduzidas com determinante igual a $79.$\\
\n  Usaremos o primeiro método.\\

\n $2\sqrt{79}\approx 17,77$, então os  possíveis valores  de  a devem ser menores que $17,77$, destes valores escolhemos os $a$ para os quais  79 é  resíduo quadrático modulo $a$. Logo, 

$a= \mp 1,2,3,5,6,7,9,10,13,14,15.$

\n Os valores de $b$ são os valores positivos da expressão $\sqrt{D}\ mod\  a$ situados entre $\sqrt{D}$ e $\sqrt{D} \mp a.$\\

\begin{itemize}
\item[•] $a=1$ , $b$ deve estar entre  $\sqrt{D} - a \approx 7,88$ e $\sqrt{79}\approx8,88$. \\ 
\n $b=8$, $c=\frac{8^2-79}{1}=-15. \Rightarrow \textcolor{blue}{(1,8,-15)}.$

\item[•] $a=2$ , $b$ deve estar entre  $\sqrt{D} - 2 \approx 6,88$ e $\sqrt{79}\approx8,88$. \\ 
\n $b=7$, $c=\frac{7^2-79}{2}=-15. \Rightarrow \textcolor{blue}{(2,7,-15)}.$

\item[•] $a=3$, $b$ deve estar entre  $\sqrt{D} - 3\approx 5,88$ e $\sqrt{79}\approx8,88$. \\ 
\n $b=7$, $c=\frac{7^2-79}{3}=-10. \Rightarrow \textcolor{blue}{(3,7,-10)}.$\\
\n $b=8$, $c=\frac{8^2-79}{3}=-5. \Rightarrow \textcolor{blue}{(3,8,-5)}.$

\item[•] $a=5$, $b$ deve estar entre  $\sqrt{D} - 5 \approx 3,88$ e $\sqrt{79}\approx 8,88$. \\ 
\n $b=7$, $c=\frac{7^2-79}{5}=-6. \Rightarrow \textcolor{blue}{(5,7,-6)}.$\\
\n $b=8$, $c=\frac{8^2-79}{5}=-3. \Rightarrow \textcolor{blue}{ (5,8,-3)}.$

\item[•] $a=6$, $b$ deve estar entre  $\sqrt{D} - 6 \approx 2,88$ e $\sqrt{79}\approx 8,88$. \\ 
\n $b=5$, $c=\frac{5^2-79}{6}=-9. \Rightarrow \textcolor{blue}{(6,5,-9)}.$\\
\n $b=7$, $c=\frac{7^2-79}{6}=-5. \Rightarrow \textcolor{blue}{ (6,7,-5)}.$

\item[•] $a=7$, $b$ deve estar entre  $\sqrt{D} - 7 \approx 1,88$ e $\sqrt{79}\approx 8,88$. \\ 
\n $b=3$, $c=\frac{3^2-79}{7}=-10. \Rightarrow \textcolor{blue}{(7,3,-10)}.$\\
\n $b=4$, $c=\frac{4^2-79}{7}=-9. \Rightarrow \textcolor{blue}{ (7,4,-9)}.$

\item[•] $a=9$, $b$ deve estar entre  $\sqrt{D} - 9 \approx -0,11$ e $\sqrt{79}\approx 8,88$. \\ 
\n $b=4$, $c=\frac{4^2-79}{9}=-7. \Rightarrow \textcolor{blue}{(9,4,-7)}.$\\
\n $b=5$, $c=\frac{5^2-79}{9}=-6. \Rightarrow \textcolor{blue}{ (9,5,-6)}.$

\item[•] $a=10$, $b$ deve estar entre  $\sqrt{D} - 10 \approx -1,11$ e $\sqrt{79}\approx 8,88$. \\ 
\n $b=3$, $c=\frac{3^2-79}{10}=-7. \Rightarrow \textcolor{blue}{(10,3,-7)}.$\\
\n $b=7$, $c=\frac{7^2-79}{10}=-3. \Rightarrow \textcolor{blue}{ (10,7,-3)}.$

\item[•] $a=13$, $b$ deve estar entre  $\sqrt{D} - 13 \approx -3,11$ e $\sqrt{79}\approx 8,88$. \\ 
\n $b=1$, $c=\frac{1^2-79}{13}=-6.$\\
Neste caso, $a=13$ não se encontra entre $\sqrt{79}-1 \approx 7,88$  e $\sqrt{79}+1\approx 9,88$. Assim, não consideramos a forma \textcolor{red}{ $(13,1,-6)$}.

\item[•] $a=14$, $b$ deve estar entre  $\sqrt{D} - 14 \approx -4,11$ e $\sqrt{79}\approx 8,88$. \\ 
\n $b=3$, $c=\frac{3^2-79}{14}=-5.$\\
Neste caso, $a=14$ não se encontra entre $\sqrt{79}-3 \approx 5,88$  e $\sqrt{79}+3\approx 11,88$. Assim, não consideramos a forma \textcolor{red}{ $(14,3,-5)$}.

\item[•] $a=15$, $b$ deve estar entre  $\sqrt{D} - 15 \approx -4,11$ e $\sqrt{79}\approx 8,88$. \\ 
\n $b=2$, $c=\frac{2^2-79}{15}=-5.$\\
Neste caso, $a=15$ não se encontra entre $\sqrt{79}-2 \approx 6,88$  e $\sqrt{79}+2\approx 10,88$. Assim, não consideramos a forma \textcolor{red}{ $(15,2,-5)$}.

\n $b=7$, $c=\frac{7^2-79}{15}=-2.$ $\Rightarrow$ \textcolor{blue}{$(15,7,-2)$}.\\
\n $b=8$, $c=\frac{8^2-79}{15}=-1.$ $\Rightarrow$ \textcolor{blue}{$(15,8,-1)$}.\\
\end{itemize}

\n Das anteriores formas resultam também:

\begin{center}
$\textcolor{blue}{(-1,8,15),(-2,7,15),(-3,7,10),(-3,7,10),(-5,7,6) ,}$\\
 $\textcolor{blue}{(-5,8,3),(-5,8,3),(-6,5,9),(-6,7,5), (-7,3,10)}$\\
 
 $\textcolor{blue}{(-7,4,9),(-9,5,6), (-10,3,7), (-10,7,3),(-15,7,2)(-15,8,1)}.$
\end{center} 

\n Portanto, temos 32 formas reduzidas com determinante $79.$
\end{exemplo}

\begin{definition}

\n Sejam $F$ uma forma reduzida com determinante $D$ e a forma $F_1$ contígua pela ultima parte a $F$; $F_2$ uma forma reduzida contígua pela ultima parte a $F_1$; $F_3$ uma  forma reduzida contígua a forma $F_2$, etc. Então, todas as formas $F_1,F_2,F_3,\cdots$ estarão completamente determinadas e serão propriamente equivalentes entre si e propriamente equivalentes  à forma $F$. Dado que o número de formas  reduzidas de um  determinante dado é finito, todas as formas da progressão  $F_1,F_2,F_3,\cdots$  não podem ser diferentes. Suponhamos que $F_m$ e $F_{m+n}$ são iguais, então $F_{m-1}$ e $F_{m+n-1}$ serão reduzidas, contíguas à mesma forma reduzida e portanto iguais. Da mesma maneira $F_{m-2}$ e $F_{m+n-2}$ serão idênticas,... e finalmente $F$ e $F_n$.  Assim, na progressão $F_1,F_2,F_3,\cdots$  necessariamente vai aparecer novamente a primeira forma $F$.\\
Suponhamos que $F_n$ é a primeira forma igual a $F$, isto é, todas as formas  $F_1,F_2,F_3, \cdots  $ $F_{n-1}$ são diferentes da forma $F$ e diferentes entre si. Chamaremos a  este conjunto de formas \textbf{o período}  de $F$.

\end{definition}
\begin{exemplo}
Determinar o período da forma $(3,8,-5)$\\
O determinante da forma é 79.\\
Seja $\textcolor{blue}{F=(3,8,-5)}$

\n Encontremos todas as formas reduzidas propriamente equivalentes a $F$

\begin{itemize}
\item[•] Seja $F_1=(a_1,b_1,a_2), a_1=-5,$\\
$ -8\equiv 7\ mod\ 5 \Rightarrow b_1= 7$ e $a_2=\frac{7^2-79}{-5}=6.$\\
Assim. $\textcolor{blue}{F_1=(-5,7,6)}$

\item[•] Seja $F_2=(a_2,b_2,a_3)$\\
$-7\equiv 5\ mod\ 6 \Rightarrow b_2= 5, a_3=-9.$\\
Assim. $\textcolor{blue}{F_2=(6,5,-9)}.$

\item[•] Seja $F_3=(a_3,b_3,a_{4})$\\
$-5\equiv 4\ mod\ 9 \Rightarrow b_3= 4, a_{4}=7.$\\
Assim. $\textcolor{blue}{F_3=(-9,4,7)}.$

\item[•] Seja $F_{4}=(a_{4},b_{4},a_{5})$\\
$-4\equiv 3\ mod\ 7 \Rightarrow b_{4}= 3, a_{5}=-10.$\\
Assim. $\textcolor{blue}{F_{4}=(7,3,-10)}.$

\item[•] Seja $F_{5}=(a_{5},b_{5},a_{6})$\\
$-3\equiv 7\ mod\ 10 \Rightarrow b_{5}= 3, a_{6}=3.$\\
Assim. $\textcolor{blue}{F_{5}=(-10,7,3)}.$

\item[•] Seja $F_{6}=(a_{6},b_{6},a_{7})$\\
$-7\equiv 8\ mod\ 3 \Rightarrow b_{6}= 8, a_{7}=-5.$\\
Assim. $\textcolor{blue}{F_{6}=(3,8,-5)=F}.$

Logo, o período de $F$ é: $(3,8,-5),(-5,7,6),(6,5,-9),(-9,4,7),(7,3,-10),$\\
$(-10,7,3).$

\end{itemize}
\end{exemplo}
\subsection{Caráter particular de uma forma quadrática primitiva.}
\begin{definition}
Uma  forma quadrática $(a,b,c)$ é \textbf{primitiva} se  os números $a,b,c$ são relativamente primos.
\end{definition}

\begin{theorem} \label{T3}
Seja $F$ uma forma  primitiva com determinante $D$ e $p$ um número primo que divide  $D$  então os números não divisíveis por $p$ que podem ser representados pela forma $F$ são todos resíduos quadráticos de $p$ ou  são todos não resíduos. 
\end{theorem}
\n \textbf{Demonstração.}

\n Sejam $m$ e $m_1$ dois números que não são  divisíveis por $p$ e  que podem  ser representados  por $F$, isto é, 

\[m= ag^2+2bgh+ch^2, \ \ \ \ \ \ \ \ \ \ m_1= ag_1^2+2bg_1h_1+ch_1^2\] 

\n Então, 
\[mm'= [agg_1+b(gh_1+hg') + chh_1]^2 -D(gh_1-hg_1)^2.\] 

\n Portanto, $mm_1$  é um resíduo quadrático módulo $D$, e portanto também é um resíduo quadrático módulo $p$, pois $p$ divide a $D$.  Daqui vem que,  $m$ e $m_1$ são ambos resíduos quadráticos, ou ambos não resíduos quadráticos módulo $p.$ $ \ \ \ \ \ \ \ \ \ \ \ \ \ \ \ \ \ \ \ \  \ \ \ \ \ \  \ \  \ \square$
\begin{observacao}
De maneira similar ao teorema anterior podemos provar o seguinte:

\begin{itemize}
\item[•] Se o determinante $D$ da forma é divisível por $4$, todos os números ímpares que podem ser representados pela forma são todos  ou $\equiv 1$ ou todos  $\equiv 3$ (mod $4$). De fato, o produto de dois desses números é resíduo quadrático de $4$ e portanto $\equiv 1\ (mod\ 4)$. Assim, ou todos são  $\equiv 1 \ mod\ 4$ ou todos são $\equiv 3 \ (mod\ 4)$.

\item[•] Se o determinante $D$ da forma é divisível por $ 8 $, o produto de qualquer dois números ímpares que podem ser representados pela forma será resíduo quadrático de $ 8 $ e portanto $\equiv 1\ (mod \ 8)$.  Assim, serão  todos congruentes a $ 1 $  ou todos congruentes a $ 3 $  ou todos congruentes a $ 5 $ ou todos congruentes a $  7$  (mod $  8$). 
\end{itemize}
\end{observacao}

\begin{theorem}
Quando o determinante $D$ da forma primitiva $F$ é $\equiv 3\ (mod\ 4)$ todos os números ímpares representáveis pela forma $F$ serão congruentes a $ 1 $ ou todos congruentes a $ 3 $ $ mod\ 4.$ 
\end{theorem}

\n \textbf{Demonstração.}

\n Sejam $m$ e $m'$ dois números que podem ser representados pela forma $F$, o produto deles pode-se reduzir à forma $p^2-Dq^2$ como vimos no Teorema \eqref{T3}. Quando cada um desses números $m,m'$ são ímpares, um dos números $p$  ou $q$  deve ser par e outro ímpar (pois $p^2-Dq^2$ é ímpar)  e portanto um dos quadrados $p^2$ e $q^2$ será congruente a $  0$ e outro congruente a $  1$ $(mod\  4). $ Assim,  $p^2-Dq^2$  deverá ser congruente a $  1$ $(mod\ 4)$, isto é, $mm' \equiv 1\ (mod\ 4)$. Logo, os dois  $m,m'$ devem ser congruentes a $ 1 $ $(mod\  4)$ ou ambos congruentes a $  3$ $(mod\ 4)$. \hspace*{10cm}$\square$
\begin{exemplo}
Consideremos a forma $(10,3,17)$ com determinante $D=-161$ o qual é $\equiv 3\ (mod\ 4)$. Dado que $17$ pode ser representado pela forma e é  $\equiv 1 \ (mod\ 4)$ então nenhum número ímpar que não seja  da forma $4n+1$ pode ser representado pela forma $(10,3,17)$.
\end{exemplo}
\begin{theorem}
Quando o determinante $D$ da forma primitiva $F$ é $\equiv 2\ (mod\ 8)$ todos os números ímpares representáveis pela forma $F$ serão  em parte  congruentes a $  1$  e em parte congruentes a $  7$, ou em parte congruentes a $  3$ e em  parte congruentes a $  5$ $(mod\ 8)$.
\end{theorem}
\begin{exemplo}
Consideremos a forma $(3,1,5)$ com determinante $D=-14$ o qual é $\equiv 2\ (mod\ 8)$. Dado que $3\ e\ 5$  são representáveis pela forma então todos os números ímpares representáveis pela forma são ou congruentes a $ 3$ ou congruentes  $  5$ módulo $8$ e nenhum número da forma $8n+1$ ou $8n+7$ pode ser representado por esta forma. 
\end{exemplo}

\begin{theorem}
Quando o determinante $D$ da forma primitiva $F$ é $\equiv 6\ (mod\ 8)$  os números ímpares representáveis pela forma $F$ são apenas aqueles que são congruentes a $1$ e $3$ $(mod\ 8)$ ou só aqueles que são  congruentes a $5$ e $7$ $(mod\ 8).$
\end{theorem}

\begin{definition} [Caráter particular de uma forma primitiva $F$]

\hspace{.2cm}

\n Quando somente resíduos quadráticos de um número primo $p$ que é  divisor  do determinante da forma $F$,  podem se representar pela forma $F$, atribuiremos  a ela o caráter  $\mathbf{Rp},$ caso contrário atribuiremos o caráter $\mathbf{Np}.$\\
\n Quando somente os números que são congruentes com $ 1 $ módulo $ 4 $ podem ser representados pela forma $F$ atribuiremos o caráter $1,4$. 
De forma similar definimos os carateres denotados por $3,4;1,8;3,8;5,8$ e $7,8.$\\
\n Se temos formas mediante as quais somente podem-se representar números ímpares congruentes  a $ 1 $ ou a $ 7 $ módulo $ 8 $, atribuiremos a elas o caráter $1\ e\ 7,8$.
\n Claramente podemos inferir o que  representam os  carateres denotados por $3\ e\ 5, 8; 1\ e\ 3,8; 5\ e\ 7,8.$
\end{definition}
\begin{observacao} \label{o1}
 \n Os diferentes carateres particulares de uma forma primitiva $(a,b,c)$ com determinante  $D$ podem-se conhecer a partir dos números $a$ e $c$.\\
\n  Com efeito,  se $p$ é um divisor primo de $D$, um dos números $a$ ou $c$ não será divisível por $p$, pois se ambos fossem divisíveis por $p$, a forma não seria primitiva pois $p$   também divide a $b$ (já que $p$ divide $b^2=D+ac$). Assim, dado que tanto $a$ como $c$ podem-se  representar mediante a forma $(a,b,c)$ então bastaria verificar se $a$ (no caso que $p$ não divida $a$) ou $c$ (no caso que $p$ não divida $c$)  é resíduo quadrático de $p$.
 
\n Nos casos em que a forma $(a,b,c)$ possui uma relação com o número $ 4 $ ou o número $ 8 $ como as descritas antes, então pelo menos um dos números $a$ ou $c$ será ímpar e usando este número ímpar poderá se saber o caráter a respeito de $4$ ou $8$. 
 
\end{observacao}
\begin{definition}
O conjunto de todos os carateres particulares de uma forma  dada constituem o \textbf{caráter completo} da forma. 
\end{definition}


\begin{exemplo} \label{ex1}
Determinar os carateres particulares da forma $(10,3,17).$

\n O determinante da forma é $D=-161$, os divisores primos de $D$ são $7$ e $23$.\\
\n  O caráter da forma $(10,3,17)$  a respeito do número $7$ pode-se deduzir do número $10$  como $N7$,  pois $10$ não é resíduo quadrático módulo    $7$. \\
\n O caráter da forma $(10,3,17)$  a respeito do número $23$ pode-se deduzir do número $10$  como $N23$,  pois $10$ não é resíduo quadrático modulo $23$. 

\n Como $D \equiv 3\ mod\ 4$, o caráter da forma $(10,3,17)$ a  respeito do número $4$ pode-se deduzir do número $17$ como $1,4$ pois  $17\equiv 1\ (mod\ 4).$

\n Logo, o caráter completo da forma $(10,3,17)$ é: $N7$;$N23$;$1,4$.
\end{exemplo}
\begin{definition}
Duas formas  dizem-se no mesmo \textbf{género} se têm o mesmo caráter completo. 
\end{definition}

\subsection{Números característicos da forma $(a,b,c).$}
\begin{definition}

\n Se a forma primitiva $F=(a,b,c)$ é tal que pode-se encontrar dois números  $g$ e $h$, tais que $g^2\equiv a, gh\equiv b, h^2\equiv c$  com respeito a um módulo $m$ dado, diremos que a \textbf{forma é um resíduo quadrático} módulo $m$ e que $gx+hy $ é o valor da expressão $\sqrt{ax^2+2bxy+cy^2}$ $mod\ m$ ou simplesmente que $(g,h)$ é o valor da expressão $\sqrt{(a,b,c)}$   ou $\sqrt{F}$ $mod\ m$. \\

\n De modo mais geral, se temos um número $M$ relativamente primo ao módulo $m$ tal que 
\[g^2 \equiv aM,\ gh\equiv bM,\ h^2\equiv cM\ \ mod\ m\]

\n dizemos que $M\cdot (a,b,c)$ ou $MF$ é um resíduo quadrático de $m$ e que $(g,h)$ é o valor da expressão  $\sqrt{M(a,b,c)}$ ou $\sqrt{MF}$ $mod\ m$. 

\end{definition}

\begin{exemplo}
Consideremos a forma $(3,1,54)$, os números $7\ e\ 10$ satisfazem 
\[7^2 \equiv 3,\ 70\equiv 1,\ 10^2\equiv 54\ \ (mod\ 23)\]
Logo, a forma é um resíduo quadrático modulo 23 e $7x+10y$ é um valor da expressão $\sqrt{(3,1,54)}.$
\end{exemplo}
\begin{proposicao}
Se $M\cdot (a,b,c)$ é resíduo quadrático modulo $m$, $m$ será um divisor do determinante da forma $(a,b,c).$ 

\end{proposicao}

\begin{proposicao}
Se $M\cdot (a,b,c)$ é resíduo quadrático modulo $m$, onde $m$ ou é um número primo ou uma potencia de um número primo, $p^\mu$, o caráter particular da forma $(a,b,c)$ com respeito do número $p$ será $Rp$ ou $Np$ segundo $M$ seja resíduo ou não resíduo de $p$.
\end{proposicao}

\n \textbf{Demonstração.}

\n Como $M \cdot (a,b,c)$ é resíduo quadrático modulo $m$, então $aM$ e $cM$ são resíduos de $m$ ou $p$. \\
Sabemos pela observação \eqref{o1} que pelo menos um dos números $a$ ou $c$ não é divisível por $p$. Suponhamos sem perda de generalidade que $p$ não divide $a$  então se $M$ é resíduo quadrático $mod\ p$,   $a$ seria resíduo quadrático $mod\ p$. Logo, o caráter particular da forma com respeito de $p$ seria $Rp$.
Agora, se $M$ é um não resíduo quadrático $mod\ p$, então $a$ é um não resíduo $mod\ p$. Logo, o caráter particular da forma com respeito de $p$ seria $Np$.
 $\square$
\begin{observacao} 
\hspace*{.1 cm}

\begin{itemize}

\item[•] Se na proposição anterior $m=4$ então ou $1,4$ ou $3,4$  será o caráter particular da forma $(a,b,c)$  segundo $M\equiv 1$ ou   $M\equiv 3\ mod\ 4$.

\item[•] Se na proposição anterior $m=8$ ou uma potencia maior do número 2 então  $1,8; 3,8;5,8;7,8$   serão os carateres  particulares da forma $(a,b,c)$  segundo $M\equiv 1;3;5;7\ mod\ 8$, respetivamente.
\end{itemize}

\end{observacao}
\begin{proposicao}\label{p3}
Se $m$ é um número primo  ou uma potencia $p^\mu$ de um número primo ímpar o qual é divisor do determinante $b^2-ac$. Se $M$ é um resíduo  ou um não resíduo de $p$ conforme o caráter particular da forma $(a,b,c)$ com respeito de $p$ seja $Rp$ ou $Np$ respetivamente, então $M(a,b,c)$ será um resíduo quadrático de $m$.
\end{proposicao}

\n \textbf{Demonstração.}

\n Suponhamos que $a$ não é divisível por $p$, então se $a$ resíduo  teremos que $aM$ é resíduo quadrático de $p$ e se se $a$ é não resíduo também  teremos que $aM$ é resíduo quadrático de $p$ e portanto será resíduo de $m$ , seja $g$ o valor da expressão $\sqrt{aM}\ mod\ m,$ $h$ um valor da expressão $\frac{bg}{a} \ (mod\ m)$, isto é,  $g^2\equiv aM\ mod\ m$, $ah\equiv bg\ (mod\ m)
 $. Então
\[agh \equiv bg^2 \equiv abM\ \ \ \ \ logo\ \ \ \ gh\equiv bM\]
e
\[ ah^2\equiv bgh\equiv b^2M \equiv b^2M - (b^2-ac)M \equiv acM\]
\n Assim, $h^2\equiv cM$. Logo, podemos concluir que $(g,h)$ é o valor da expressão $\sqrt{M(a,b,c)}$ $(mod\ m)$, isto é, $M(a,b,c)$ é resíduo quadrático $(mod\ m)$. Se $a$ for divisível por $p$ então certamente $c$ não será divisível por $p$ e procederíamos da mesma forma anterior e obtendo o mesmo resultado.$\hspace{9.3cm}\square$
\begin{observacao}
De forma similar a proposição anterior pode-se demonstrar que  se $m=4$ e $m$ divide  $b^2-ac$, e se o numero $M$ se considera ou $\equiv 1$ ou $\equiv 3$ segundo ou $1,4$ ou $3,4$ seja o caráter particular da forma $(ab,c) $ então $M(a,b,c)$ é resíduo quadrático de $m$. Se $m=8$ ou uma potencia maior de 2 e $m$ divide a $b^2-ac$, e se $M\equiv 1;3;5;7\ mod\ 8$ segundo seja o caráter da forma respeito do número $8$,  então  $M(a,b,c)$ é resíduo quadrático de $m$.
\end{observacao}

\n Se o determinante da forma $(a,b,c)$ é igual a $D$ e $M(a,b,c)$ é resíduo quadrático de $D$, a partir do número $M$ pode-se encontrar todos os carateres particulares da forma $(a,b,c)$ em respeito a todos os divisores primos ímpares de $D$ em respeito ao número 4 e 8 (se eles dividem a $D$).\\

\n Por exemplo, consideremos a forma $(20,10,27)$ com determinante $-440$ $(440=2^3\cdot 5 \cdot11).$
Dado que, $3(20,10,27)$ é resíduo quadrático de 440, onde  $(150,9)$ é o valor da expressão $\sqrt{3(20,10,27)}$  e 3 é um não resíduo módulo 5 e 3 é um resíduo de 11, então os carateres particulares da forma $(20,10,27)$ em  respeito a 5 e  a 11 são $N5$ e $R11$, respetivamente e em respeito a 8 é $3,8$.

\begin{observacao}
Os carateres particulares em  respeito a $4$ e a $8$ , sempre que eles não sejam divisores do determinante, são os únicos que não necessariamente  têm uma conexão  com o número $M$.
\end{observacao}

\begin{proposicao}
Se o número $M$ é primo relativo com o determinante $D$  e de $M$ podemos deduzir todos os carateres particulares da forma $(a,b,c)$ (exceto aqueles carateres respeito de $  4$ e $  8$, onde $  4$ e $  8$ não dividem $D$)  então $M(a,b,c)$ será um resíduo quadrático de $D$. 
\end{proposicao}
\n \textbf{Demonstração.}

\n Se escrevemos $D$ da forma $\pm p_1^{k_1} p_2^{k_2} \cdots p_s^{k_s}  $ onde $p_1,p_2,\cdots p_s$ são primos distintos, pela proposição \ref{p3}, $M(a,b,c)$ será um resíduo quadrático de cada um dos $p_1^\alpha,p_2^\beta, p_3^\gamma,  \cdots  $\\
Suponha que o valor de $\sqrt{M(a,b,c)}$ é $(u_1,v_1)$ em respeito a $mod\ p_1^{k_1}$, em respeito a $mod\ p_2^{k_2}$ é $(u_2,v_2)$,$ \cdots$ , em respeito a $mod\ p_s^{k_s}$ é $(u_s,v_s)$.\\
Se os números $g,h$ são tais que $g \equiv u_1,u_2,\cdots, u_{s}$ e $h\equiv v_1,v_2,\cdots, v_s$ em respeito dos módulos $p_1^{k_1},p_2^{k_2}, \cdots, p_s^{k_s}$, respetivamente. Então\\

\begin{multicols}{2}
$\begin{array}{llllllll}
g^2&\equiv & u_1^2& \equiv & aM & mod\ p_1^{k_1}\\
g^2&\equiv & u_2^2& \equiv & aM& mod\ p_2^{k_2}\\
&&&\vdots&& \\
g^2&\equiv & u_s^2 &\equiv & aM& mod\ p_s^{k_s}\\

\end{array}$ \\ $\begin{array}{llllllll}
gh&\equiv & u_1v_1& \equiv & bM & mod\ p_1^{k_1}\\
gh&\equiv & u_2v_2& \equiv & bM& mod\ p_2^{k_2}\\
&&&\vdots&& \\
gh&\equiv & u_sv_s&\equiv & bM& mod\ p_3^{k_s}\\

\end{array}$
\end{multicols}

\[\begin{array}{llllllll}
h^2&\equiv & v_1^2& \equiv & cM & mod\ p_1^{k_1}\\
h^2&\equiv & v_2^2& \equiv & cM& mod\ p_2^{k_2}\\
&&&\vdots&& \\
h^2&\equiv & v_s'^2 &\equiv & cM& mod\ p_3^{k_s}\\

\end{array}\]

e portanto, 

\begin{multicols}{2}
$\begin{array}{llllllll}
g^2&\equiv & aM & mod\  p_1^{k_1} p_2^{k_2} \cdots p_s^{k_s}     \\
gh&\equiv& bM & mod\  p_1^{k_1} p_2^{k_2} \cdots p_s^{k_s}  & \Leftrightarrow \\
h^2&\equiv & cM & mod\  p_1^{k_1} p_2^{k_2} \cdots p_s^{k_s}   \\

\end{array}$\\ 
$\begin{array}{llllllll}
g^2&\equiv & aM & mod\ D\\
gh&\equiv& bM & mod\ D  \\
h^2& \equiv & cM & mod\  D\\ 

\end{array}$
\end{multicols}

\n Assim, $M(a.b.c)$ é resíduo quadrático de $D$.$\ \ \ \ \ \  \ \ \ \ \ \ \ \ \ \ \ \ \ \ \ \ \ \ \ \ \ \ \ \ \ \ \ \ \ \ \ \ \ \ \ \ \ \ \ \ \ \  \square$\\

\n Números como $M$ são chamados de \textbf{números característicos da forma $(a,b,c).$}
\begin{observacao}
\hspace*{.1cm}
\begin{itemize}
  \item[•] Se $M$ é um número característico de uma forma primitiva com determinante $D$ dado, todos os números congruentes a $M$ modulo $D$ serão números característicos da mesma forma.  
\item[•] Formas do mesmo género possuem os mesmos números característicos.   
 
\end{itemize}

\end{observacao}
\subsection{Composição de formas}

\begin{definition}
Seja $F=AX^2+2BXY+CY^2$  uma forma que pode-se transformar no produto de duas formas $f_1=a_1x_1^2+2b_1x_1y_1+c_1y_1^2$ e $f_2=a_2x_2^2+2b_2x_2y_2+c_2y_2^2$ através da  substituição:
\[X=p_1x_1x_2+p_2x_1y_2+p_3y_1x_2+p_4y_1y_2\].
\[Y=q_1x_1x_2+q_2x_1y_2+q_3y_1x_2+q_4y_1y_2\].
\n Então, $F$ diz-se \textbf{transformável}  em $ff'$. Se além disso, os números 
\[p_1q_2-q_1p_2, p_1q_3-q_1p_3,p_1q_4-q_1p_4,p_2q_3-q_2p_3,p_2q_4-q_2p_4,p_3q_4-q_3p_4,\]
\n não possuem um divisor comum, chamaremos a forma $F$ de forma \textbf{composta} das formas $f_1$ e $f_2$.
\end{definition}
\n Para simplificar a anterior definição, diremos que  $F$  transforma-se em $f_1f_2$ através da substituição $p_1,p_2,p_3,p_4,q_1,q_2,q_3,q_4$.\\

\n A continuação apresentaremos algumas conclusões que se obtém a partir  da suposição de que a forma $F$ se transforma em $f_1f_2$ através da substituição $p_1,p_2,p_3,p_4,q_1,q_2,q_3,q_4$. Estas conclusões apresentam-se com detalhe em \cite{1}.

\n Antes definiremos alguns elementos. \\
\n Sejam $D,d_1\ d_2$  os determinantes das formas $F,f_1,f_2$ respetivamente e $M,m_1,m_2$ os máximos divisores comuns dos números $A,2B,C;a_1,2b_1,c_1; a_2,2b_2,c_2$, respetivamente. \\
\n Além disso, sejam $\alpha_1,\beta_1,\gamma_1,\alpha_2,\beta_2,\gamma_2$ inteiros tais que:
\[\alpha_1 a+2\beta_1 b+\gamma_1c=m_1, \ \ \ \ \ \ \ \ \ \ \ \ \ \alpha_2a_2+2\beta_2b_2 +\gamma_2 c_2=m_2.\]
\n Denotemos os números 
\[p_1q_2-q_1p_2, p_1q_3-q_1p_3,p_1q_4-q_1p_4,p_2q_3-q_2p_3,p_2q_4-q_2p_4,p_3q_4-q_3p_4,\]
por $P,Q,R,S,T,U$ respetivamente e seja $k$ o máximo divisor comum deles  tomado positivamente. 

\begin{itemize}

\item[1.] Os  determinantes das formas $F,f_1,f_2$ satisfazem $d_1=Dn_1^2, d_2=Dn_2^2$ \footnote{$\alpha_1 P + \beta_1 (R-S)+  \gamma_1 U=m_1n_2, \alpha_2  Q + \beta_2 (R+S)+  \gamma_2T=m_2n_1,  $ onde $n_1 \ e \ n_2$ podem ser frações sempre que $m_1n_2$ e $m_2n_1$ sejam inteiros.} )
\item[2.]  $D$ sempre divide os números $d_1m_2^2\ e\  d_2m_1^2$. Logo, $D,d_1,d_2$ tem o mesmo sinal e não existe uma forma que se possa transformar no produto $f_1f_2$ se seu determinante é maior que o máximo divisor comum dos números $d_1m_2^2$ e $d_2m_1^2$.

\item[3.] $k$ é o máximo divisor comum dos números $m_1n_2$ e $m_2n_1$ e $Dk^2$ será o máximo divisor comum dos números $d_1m_2^2$ e $d_2m_1^2$.

\item[4.] Se $F$ é composta pelas formas $f_1$ e $f_2$, isto é $k=1$, $M$ necessariamente será igual a $m_1m_2$.

\end{itemize}
\n \textbf{Problema:} Dadas duas formas $f_1,f_2$ cujos determinantes sejam iguais ou  que no máximo diferem por fatores quadrados encontrar uma forma composta por estas duas formas. \\
\n Sejam as formas $f_1=(a_1,b_1,c_1)$ e $f_2=(a_2,b_2,c_2)$ com determinantes $d_1,d_2$ respetivamente e sejam $m_1$ e $m_2$ os máximos divisores comuns dos números $a_1,2b_1,c_1; a_2,2b_2 ,$ $c_2$ respetivamente.\\
\n Seja  $D$ o máximo divisor comum dos números  $d_1m_2^2$ e $d_2m_1^2$  tomados com o mesmo sinal que $d_1$ e $d_2$.  Então $\frac{d_1m_2^2}{D}$  e $\frac{d_2m_1^2}{D} $ serão números positivos relativamente primos cujo produto será um quadrado e portanto cada um deles será um quadrado. Logo, $\sqrt{\frac{d_1}{D}}$ e $\sqrt{\frac{d_2}{D} } $ são números  racionais os quais denotaremos por  $n_1,n_2$ respetivamente. Os números $m_1n_2$ e $m_2n_1$ são relativamente primos e os números $a_1n_2,c_1n_2,a_2n_1,b_1n_2+b_2n_1$ e $b_1n_2-b_1n_1$ são inteiros.  \\
\n Agora consideremos quatro inteiros $\mathcal{Q} _1,\mathcal{Q}_2,\mathcal{Q}_3$ e $\mathcal{Q}_4$ arbitrários  com uma condição, que as quatro quantidades do lado esquerdo da seguinte equação \eqref{eI} não sejam todos iguais a zero. Sejam as equações

\begin{equation}\label{eI}
\begin{array}{rcl}
\mathcal{Q}_2a_1n_2+\mathcal{Q}_3a_2n_1+\mathcal{Q}_4(b_1n_2+b_2n_1)&=& \mu q_1\\
-\mathcal{Q}_1a_1n_2+\mathcal{Q}_4c_2n-\mathcal{Q}_3(b_1n_2-b_2n_1)&=& \mu q_2\\
\mathcal{Q}_4c_1n_2-\mathcal{Q}_1a_2n_1+\mathcal{Q}_2(b_1n_2-b_2n_1)&=& \mu q_3\\
-\mathcal{Q}_3c_1n_2-\mathcal{Q}_2c_2n_1-\mathcal{Q}_1(b_1n_2+b_2n_1)&=& \mu q_4\\
\end{array}
\end{equation} 

\n tais que $q_1,q_2,q_3\ e \ q_4$ são inteiros que não possuem divisores comuns. Isto pode-se lograr tomando  $\mu$  como o máximo comum divisor dos quatro números que estão a esquerda das equações. Como $mdc\ (q_1,q_2,q_3,q_4)=1$ então podemos encontrar inteiros $\mathcal{B}_1,\mathcal{B}_2,\mathcal{B}_3,\mathcal{B}_4$ tais que 
\[\mathcal{B}_1 q_1+\mathcal{B}_2q_2+\mathcal{B}_3q_3+\mathcal{B}_4 q_4=1\].
\n Tendo isto, determinamos os números $p_1,p_2,p_3,p_4$ mediante as seguintes equações:

\begin{equation}\label{eII}
\begin{array}{rcl}
\mathcal{B}_1a_1n_2+\mathcal{B}_3 a_2n_1 + \mathcal{B}_4(b_1n_2+b_2n_1)&=&p_1\\
-\mathcal{B}_1 a_1n_2+\mathcal{B}_4 c_2n_1 - \mathcal{B}_3(b_1n_2-b_2n_1)&=&p_2\\
\mathcal{B}_4c_1n_2-\mathcal{B}_1 a_2n_1 + \mathcal{B}_2(b_1n_2-b_2n_1)&=&p_3\\
-\mathcal{B}_3c_1n_2-\mathcal{B}_2 c_2n_2 - \mathcal{B}_4(b_1n_2+b_2n_1)&=&p_4\\

\end{array}
\end{equation} 

\n Agora, fazemos
\[q_2q_3-q_1q_4=An_2n_2, \ \ p_1q_4+q_1p_4-p_2q_3-q_2p_3=2Bn_1n_2,\ \ p_2p_3-p_1p_4=Cn_1n_2.\]

\n Então, os números $A,B$ e $C$ são inteiros e a forma $F=(A,B,C)$ será a forma composta pelas formas $f_1$ e $f_2$.

\n Agora, vamos a limitar o problema anterior. Consideraremos as formas $f_1$ e $f_2$ com o mesmo determinante, isto é; $d_1=d_2$, $m_1$ e $m_2$ devem ser relativamente primos e a forma $F$ deve ser composta diretamente por $f_1$ e $f_2$. \\
\n Como $(m_1,m_2)= 1$ então $m_1^2\ e\ m_2^2$ são também relativamente primos e assim $D$ que é o máximo divisor comum dos números $d_1m_2^2$ e $d_2m_1^2$, , será igual a $d_1(=d_2)$  e $n_1=n_2=1$. \\
\n Como os valores de $\mathcal{Q}_1,\mathcal{Q}_2,\mathcal{Q}_3,\mathcal{Q}_4$ são arbitrários então consideremos  eles iguais a $-1,0,0,0$ respetivamente. Logo, as equações \eqref{eI} ficam assim:

\begin{equation*}
\begin{array}{rcl}
0 &=& \mu q_1\\
a_1&=& \mu q_2\\
a_"&=& \mu q_3\\
b_1+b_2&=& \mu q_4\\
\end{array}
\end{equation*} 

\n Logo, 

\[\mathcal{B}_1 \cdot 0+\mathcal{B}_2 q_2 + \mathcal{B}_3q_3+\mathcal{B}_4q_4=1\]
\[\mathcal{B}_2 q_1 \mu+ \mathcal{B}_3 q_3\mu +\mathcal{B}_4q_4\mu =\mu\]
\[\mathcal{B}_2 a+ \mathcal{B}_3 a_2 +\mathcal{B}_4(b_1+b_2)=\mu\]

 \n Agora, se  $\mu$ é o máximo comum divisor dos números $a_1,a_2,b_1+b_2$ então os números $\mathcal{B}_2,\mathcal{B}_3,\mathcal{B}_4$ se podem escolher de forma que:
 \[\mathcal{B}_2 a_1 +\mathcal{B}_3 a_2+\mathcal{B}_4(b_1+b_2)=\mu\]
\n e $\mathcal{B}_1$ pode-se escolher arbitrariamente. \\

\n As equações \eqref{eII} ficam a ser

\begin{equation*}
\begin{array}{rcl}
\mathcal{B}_2a_1+\mathcal{B}_3 a_2+\mathcal{B}_4(b_1+b_2)&=&p_1\\
-\mathcal{B}_1 a_1+\mathcal{B}_4 c_2- \mathcal{B}_3(b_1-b_2)&=&p_2\\
\mathcal{B}_4c_!-\mathcal{B}_1 a_2 + \mathcal{B}_2(b_1-b_2)&=&p_3\\
-\mathcal{B}_3c_1-\mathcal{B}_2 c_2 - \mathcal{B}_4(b_1+b_2)&=&p_4\\

\end{array}
\end{equation*} 

\n Os valores de A e B são:
\[A= q_2q_3 =\dfrac{a_1a_2}{\mu^2}\]
\begin{equation*}
\begin{array}{lll}
2B&=&p_1q_4+q_1p_4-p_2q_3-q_2p_3\\
2B&=& \frac{1}{\mu} \left[ p_1(b_1+b_2)-p_2a_2-p_3a_1\right] \\
2B&=& \frac{1}{\mu} \left[ (\mathcal{B}_2a_1+\mathcal{B}_3 a_2+ \mathcal{B}_4(b_1+b_2))(b_1+b_2)-(-\mathcal{B}_1 a_1+\mathcal{B}_4 c_2- \mathcal{B}_3(b_1-b_2))a_2 \right. \\
&&\left. -(\mathcal{B}_4c_1-\mathcal{B}_1 a_1 + \mathcal{B}_2(b_1-b_2))a_1\right] \\
2B&=& \frac{1}{\mu} \left[2\mathcal{B}_1 a_1a_2+2\mathcal{B}_2a_1b_2+2\mathcal{B}_3a_2b_1+2\mathcal{B}_4b_1b_2+\mathcal{B}_4(b_1^2-a_1c_1) +\mathcal{B}_4(b_2^2-a_2c_2) \right] \\
2B&=& \frac{1}{\mu} \left[2\mathcal{B}_1 a_1a_2+2\mathcal{B}_2a_1b_2+2\mathcal{B}_3a_2b_1+2\mathcal{B}_4b_1b_2+\mathcal{B}_4D +\mathcal{B}_4D\right] \\
2B&=& \frac{1}{\mu} \left[2\mathcal{B}_1 a_1a_2+2\mathcal{B}_2a_1b_2+2\mathcal{B}_3a_2b_1+2\mathcal{B}_4b_1b_2+2\mathcal{B}_4D \right] \\

2B&=& \frac{1}{\mu} \left[2\mathcal{B}_1 a_1a_2+2\mathcal{B}_2a_1b_2+2\mathcal{B}_3a_2b_1+2\mathcal{B}_4(b_1b_2+D) \right] \\
&&\\
B&=& \frac{1}{\mu} \left[\mathcal{B}_1 a_1a_2+\mathcal{B}_2a_1b_2+\mathcal{B}_3a_2b_1+\mathcal{B}_4(b_1b_2+D) \right] \\
\end{array}
\end{equation*}

\n $C$ pode determinar-se mediante a equação $AC=B^2-D$ se $a_1$ e $a_2$ não são simultaneamente zero. \\

\n Do anterior podemos observar que o valor de $A$ não depende  dos valores de $\mathcal{B}_1,\mathcal{B}_2,\mathcal{B}_3,$ $\mathcal{B}_4$ mas o valor de $B$ varia se variam os valores dos $\mathcal{B}_1,\mathcal{B}_2,\mathcal{B}_3,\mathcal{B}_4$.  Pode-se mostrar que não importa como se determinem os valores de $\mathcal{B}_1,\mathcal{B}_2,\mathcal{B}_3,\mathcal{B}_4$, todos os valores de $B$ serão congruentes modulo $A$.\\

\n Da expressão
\[\mathcal{B}_2a_1 + \mathcal{B}_3 a_2+\mathcal{B}_4 (b_1+b_2)=\mu\]

\n podemos deduzir que
\[\mathcal{B}_2a_1b_1 +\mathcal{B}_3 a_2b_1+\mathcal{B}_4(b_1+b_2)b_1=b_1\mu\]
Assim,

\n $\begin{array}{llll}

B&=&\dfrac{1}{\mu} \left[\mathcal{B}_1 a_1a_2+\mathcal{B}_2a_1b_2+\mathcal{B}_3a_2b_1+\mathcal{B}_4(b_1b_2+D) \right]\\
&=&\dfrac{a_1}{\mu} \left[\mathcal{B}_1a_2+\mathcal{B}_2b_2 \right] +\frac{1}{\mu}\mathcal{B}_4(b_1b_2+D)+\dfrac{1}{\mu}\left[ \mathcal{B}_3a_2b_1 +\mathcal{B}_2a_1b_1+\mathcal{B}_4(b_1+b2)b_1\right. \\
&&\ \ \  \left. -\mathcal{B}_2a_1b_1-\mathcal{B}_4(b_1+b_2)b_1 \right] \\
&=&\dfrac{a_1}{\mu} \left[\mathcal{B}_1 a_2+\mathcal{B}_2b_2 \right] +\frac{1}{\mu}\mathcal{B}_4(b_1b_2+D)+\dfrac{1}{\mu}\left[b_1\mu- \mathcal{B}_2a_1b_1-\mathcal{B}_4(b_1+b_2)b_1\right] \\

&=&b_1+\dfrac{a_1}{\mu} \left[ \mathcal{B}_1 a_2 +\mathcal{B}_2(b_2-b_1)-\mathcal{B}_4 c_1\right] .
\end{array}$

\n Similarmente  podemos obter 

\[B= b_2+\frac{a_2}{\mu} \left[ \mathcal{B}_1 a_1+\mathcal{B}_3(b_1-b_2)-\mathcal{B}_4 c_2\right]. \]

\n Portanto, 

\[B \equiv b_1 \ mod\ \dfrac{a_1}{\mu} \ \ \  e \ \ \ \ B \equiv b_2\ mod\ \dfrac{a_2}{\mu}.  \]

\n Sempre que $\dfrac{a_1}{\mu}$, $\dfrac{a_2}{\mu}$ sejam relativamente primos podemos encontrar entre 0 e $A-1$ (ou entre 0 e $-A-1$, quando $A$ é negativo) um único número que será congruente como $b_1$ $mod \ \dfrac{a_1}{\mu} \ e\ com \ b_2\ mod\ \dfrac{a_2}{\mu}$. Se denotamos este número por $B$  e $\dfrac{B^2-D}{A}=C$ então a forma $(A,B,C) $ estará composta pelas formas $(a_1,b_1,c_1)\ e\ (a_2,b_2,c_2)$, neste caso não é necessário achar os $\mathcal{B}_1,\mathcal{B}_2,\mathcal{B}_3,\mathcal{B}_4$ para calcular a forma composta. \\

\n Por exemplo sejam as formas $(10,3,11)$ e $(15,2,7)$, calculemos a forma compostas por elas. \\

\n$a_1=10, a_2=15\ e\ b_1+b_2= 5$, $\mu =mdc(10,15,5)=5$.\\

\n $A= \frac{150}{25}=6$  e como $(10/5,15/5)=(2,3)=1$ então existirá entre $0\  e\  5$  um único número $B$ tal que $B\equiv b_1\ mod\  2$ e $B\equiv b_2\ mod\  3$, isto é, $B=5$, o valor de $C= \frac{126}{6}=21$. Portanto, a forma composta é $(6,5,21).$
\n Em seguida descreveremos um caso particular que será usada no exemplo  da factorização de um número. \\

\n  \textbf{Composição de formas cujos termos iniciais são potencias  de números primos.}\footnote{um número primo pode-se considerar como sua própria primeira potencia.} \\

\n Sejam duas formas primitivas $(a_1,b_1,c_1)$ e $(a_2,b_2,c_2)$ onde $a_1$ e $a_2$ são potencias do mesmo primo, $h$. Sejam $a_1=h^\chi$  e $a_2=h^\lambda$, vamos supor que $\chi $ é maior ou igual a $\lambda$. Assim, $h^\lambda$ será o máximo divisor comum de $a_1$ e $a_2$. Se alem disso, $h^\lambda$ é divisor de $b_1+b_2$ então $A=n^{\chi-\lambda}$, $B \equiv b_1\ mod\ h^{\chi-\lambda}$ e $B \equiv b_2\ mod\ 1$ e $C=\frac{B^2-D}{A}.$ Se $h^\lambda$ não divide $b+b'$, o máximo comum divisor de estes números será necessariamente uma potencia de $h$, digamos $h^\nu$ com $\nu < \lambda$ (onde $\nu =0$ se $h^{\lambda}$ e $b_1+b_2$ são relativamente primos). Se $\mathcal{B}_1,\mathcal{B}_2,\mathcal{B}_3,\mathcal{B}_4$ se determinam de forma que 

\[\mathcal{B}_2h^{\chi} + \mathcal{B}_3 h^{\chi} +\mathcal{B}_4(b_1+b_2)=h^\nu\]

\n com $\\mathcal{B}_1$ arbitrário, a forma $(A,B,C) $  será a forma composta pelas formas dadas se 

\[A=h^{\chi+\lambda-2\nu},\ \ B=b_1+h^{\chi-\nu}(\mathcal{B}_1 h^\lambda-\mathcal{B}_2(b_1-b_2)-\mathcal{B}_4c),\ \ C=\dfrac{B^2-D}{A}.\] 

\n Neste caso, $\mathcal{B}_2$ também pode-se escolher arbitrariamente, então fazendo $\mathcal{B}_1=\mathcal{B}_2=0$ temos

\[B=b_1-\mathcal{B}_4c_1h^{\chi-\nu}\]

\n o mais geralmente, 

\[B=kA+b_1-\mathcal{B}_4ch^{\chi-\nu}\]

\n onde $k$ é um número arbitrário.  O valor de $\mathcal{B}_4$ é o valor da expressão $\frac{h^\nu}{b_1+b_2}\ mod\ h^\lambda  $

\section{Fatorização de um número usando formas quadráticas. }
\n Sejam $(A,B,C)$ e $(A_1,B_1,C_1)$ formas quadráticas com o mesmo determinante $M$ ou $-M$ ou mais geralmente $\pm kM$ pertencentes ao mesmo género, então os números $AA_1,AC_1\ e\ A_1C$ são resíduos quadráticos de $kM$, pois qualquer número característico, digamos $m$, de uma forma  é também  número característico da outra, logo $mA,mC,mA_1\ e\ mC_1$ são resíduos quadráticos de $kM$. Então, o produto deles também  é resíduo quadrático. Assim,  $AA_1,AC_1\ e\ A_1C$ são resíduos quadráticos de $kM$.

\n Se $(a,b,a_1) $ é uma forma reduzida com determinante positivo $M$ ou mais geralmente  $kM$ e  $(a_1,b_1,a_2),(a_2,b_2,a_3), \cdots$ são formas no seu período então elas serão equivalentes a $(a,b,a_1)$  e certamente estarão contidas no mesmo género. Os números $aa_1,aa_2,aa_3,\cdots$ serão resíduos quadráticos módulo  $M$.\\

\n Assim, se queremos fatorizar um número $M$, a ideia é ter uma forma com determinante $M$ (ou $kM$), calcular um número considerável de formas  no seu período e obter resíduos quadráticos  modulo $M$.
\begin{exemplo} \label{e2.1}
\ Seja $M=997331$\\

\n Usemos  as formas quadráticas para obter resíduos quadráticos módulo $M.$
\begin{itemize}
\item[1.]  Consideremos a forma $(1,998,-1327)$   cujo determinante é $D=998^2+1327=997331$ e $\sqrt{D}\approx 998,66.$\\ 

\n Calculemos algumas formas pertencentes ao período de \\
 \[(a,b,a_1)=(1,998,-1327).\]
\begin{itemize}
\item[•] $(a,b,a_1)=\textcolor{blue}{(1,998,-1327)}$\\
$-998 \equiv 329 \ mod\ 1327 \Rightarrow b_1= 329$ \\
$a_2= \frac{329^2-997331}{-1327}=670.$\\
Logo, $(a_1,b_1,a_2)=\textcolor{blue}{(-1327,329,670)}$\\

\item[•] $(a_1,b_1,a_2)=(-1327,329,670)$\\
$-329 \equiv 341 \ mod\ 670 \Rightarrow b_2= 341$ \\
$a_3= \frac{341^2-997331}{670}=-1315.$\\
Logo, $(a_2,b_2,a_3)=\textcolor{blue}{(670,341,-1315).}$\\

\item[•] $(a_2,b_2,a_3)=(670,341,-1315)$\\
$-341 \equiv 974 \ mod\ 1315 \Rightarrow b_3= 974$ \\
$a_{4}= \frac{974^2-997331}{-1315}=37.$\\
Logo, $(a_3,b_3,a_{4})=\textcolor{blue}{(-1315,974,37).}$ \\

\item[•] $(a_3,b_3,a_{4})=(-1315,974,37)$\\
$-974 \equiv 987 \ mod\ 37 \Rightarrow b_{4}= 987$ \\
$a_{5}= \frac{987^2-997331}{37}=-626$\\
Logo, $(a_{4},b_{4},a_{5})=\textcolor{blue}{(37,987,-626).}$ \\

\item[•] $(a_{4},b_{4},a_{v})=(37,987,-626)$\\
$-987 \equiv 891 \ mod\ 626\Rightarrow b_{5}= 891$ \\
$a_{6}= \frac{626^2-997331}{-626}=325$\\
Logo, $(a_{5},b_{5},a_{6})=\textcolor{blue}{(-626,891,325).}$ \\

\item[•] $(a_{5},b_{5},a_{6})=(-626,891,325)$\\
$-891 \equiv 734 \ mod\ 891\Rightarrow b_{6}= 734$ \\
$a_{7}= \frac{734^2-997331}{325}=-1411$\\
Logo, $(a_{6},b_{6},a_{7})=\textcolor{blue}{(325,734,-1411).}$ \\

\item[•] $(a_{6},b_{6},a_{7})=(325,734,-1411)$\\
$-734\equiv 677 \ mod\ 1411\Rightarrow b_{7}= 677$ \\
$a_{8}= \frac{677^2-997331}{-1411}=382$\\
Logo, $(a_{7},b_{7},a_{8})=\textcolor{blue}{(-1411,677,382).}$ \\

\item[•] $(a_{7},b_{7},a_{8})=(-1411,677,382)$\\
$-677\equiv851 \ mod\ 382\Rightarrow b_{8}= 851$ \\
$a_{8}= \frac{851^2-997331}{382}=-715$\\
Logo, $(a_{8},b_{8},a_{9})=\textcolor{blue}{(382,851,-715).}$ \\
\end{itemize}
\n Em resumo, temos as seguintes formas que pertencem ao período de $(1,998,-1327)$.\\

$\begin{array}{rrrr}
(&1,&998,&-1327)  \\
(&-1327,&329,&670)\\
(&670,&341,&-1315)\\
(&-1315,&974,&37)\\
(&37,&987,&-626)\\
(&-626,&891,&325)\\
(&325,&734,&-1411)\\
(&-1411,&677,&382)\\
(&382,&851,&-715)\\

\end{array}$

\n Assim, destas formas obtemos os seguintes resíduos quadráticos modulo $M$.
\[-1327,670,-1315,37,-626,325,-1411,382,-715.\]

\item[2.] Consideremos a forma $(1,1412,-918)$   cujo determinante é $D=1994662=2\cdot 997331=2\cdot M$ e $\sqrt{D}\approx 1412,32.$\\ 

\n Calculemos algumas formas pertencentes ao período de  $(a,b,a_1)=(1,1412,-918).$
\begin{itemize}
\item[•]$(a,b,a_1)=\textcolor{blue}{(1,1412,-918)}$\\
$-1412\equiv 1342 \ mod\ 918 \Rightarrow b_1= 1342$ \\
$a_2= \frac{1342^2-1994662}{-918}=211.$\\
Logo, $(a_1,b_1,a_2)=\textcolor{blue}{(-918,1342,211)}$\\

\item[•] $(a_1,b_1,a_2)=(-918,1342,211)$\\
$-1342 \equiv 1401 \ mod\ 211 \Rightarrow b_2= 1401$ \\
$a_3= \frac{1401^2-1994662}{211}=-151.$\\
Logo, $(a_2,b_2,a_3)=\textcolor{blue}{(211,1401,-151).}$\\

\item[•] $(a_2,b_2,a_3)=(211,1401,-151)$\\
$-1401 \equiv 1317 \ mod\ 151 \Rightarrow b_3=  1317$ \\
$a_{4}= \frac{1317^2-1994662}{-151}=1723.$\\
Logo, $(a_3,b_3,a_{4})=\textcolor{blue}{(-151,1317,1723).}$ \\

\item[•] $(a_3,b_3,a_{4})=(-151,1317,1723)$\\
$-1317 \equiv 406 \ mod\ 1723 \Rightarrow b_{4}= 406$ \\
$a_{5}= \frac{406^2-1994662}{1723}=-1062$\\
Logo, $(a_{4},b_{4},a_{5})=\textcolor{blue}{(1723,406,-1062).}$ \\

\item[•] $(a_{4},b_{4},a_{5})=(1723,406,-1062)$\\
$-406 \equiv 656 \ mod\ 1062\Rightarrow b_{5}= 656$ \\
$a_{6}= \frac{656^2-1994662}{-1062}=1473$\\
Logo, $(a_{5},b_{5},a_{6})=\textcolor{blue}{(-1062,656,1473).}$ \\

\item[•] $(a_{5},b_{5},a_{6})=(-1062,656,1473)$\\
$-656\equiv 817 \ mod\ 1473\Rightarrow b_{6}= 817$ \\
$a_{6}= \frac{817^2-1994662}{1473}=-901$\\
Logo, $(a_{6},b_{6},a_{7})=\textcolor{blue}{(1473,817,-901).}$ \\

\item[•] $(a_{6},b_{6},a_{7})=(1473,817,-901)$\\
$-817\equiv 985 \ mod\ 901\Rightarrow b_{7}= 985$ \\
$a_{8}= \frac{985^2-1994662}{-901}=1137$\\
Logo, $(a_{7},b_{7},a_{8})=\textcolor{blue}{(-901,985,1137).}$ 

\end{itemize}
\end{itemize}
\n Portanto, temos as seguintes formas que pertencem ao período de $(1,1412,-918)$.\\

$\begin{array}{rrrr}
(&1,&1412,&-918)\\
(&-918,&1342,&211)\\
(&211,&1401,&-151)\\
(&-151,&1317,&1723)\\
(&1723,&406,&-1062)\\
(&-1062,&656,&1473)\\
(&1473,&817,&-901)\\
(&-901,&985,&1137)
\end{array}$\\

\n Assim, das formas obtemos os seguintes resíduos quadráticos modulo $M$.
\[-918,211,-151,1723,-1062,1473,-901,1137.\]
\end{exemplo}

\n Dos resíduos encontrados na parte $(1)$ e $ (2)$ do exemplo, consideremos os que tem fatores primos não muito grandes e usemos o fato que o produto de dois resíduos quadráticos gera outro resíduo quadrático para encontrar novos resíduos quadráticos.\\

\begin{itemize}
\item[•] $\mathbf{\textcolor{OliveGreen}{37}}.$
\item[•]$670=\mathbf{\textcolor{OliveGreen}{2\cdot 5\cdot67}}.$
\item[•] $-1411=\mathbf{\textcolor{OliveGreen}{-17\cdot83}}.$
\item[•] $-715= \mathbf{\textcolor{OliveGreen}{-5\cdot11\cdot13}}.$
\item[•] $-901=\mathbf{\textcolor{OliveGreen}{-17\cdot 53}}.$
\item[•] Como $-1062=-2\cdot 3^2\cdot59$ é resíduo quadrático de $M$ então 

$\left( \frac{-1062}{M}\right) =1 \Rightarrow \left( \frac{-2}{M}\right) \left( \frac{3^2}{M}\right) \left( \frac{59}{M}\right) =1 \Rightarrow \left( \frac{-2}{M}\right)\left(  \frac{59}{M}\right) =1 $

Logo, $\mathbf{\textcolor{OliveGreen}{-2\cdot 59}}$ é resíduo quadrático. 

\item[•] Como $-918=-2\cdot 3^3\cdot17$  resíduo quadrático de $M$ então  $\mathbf{\textcolor{OliveGreen}{-2\cdot3\cdot 17}}$ é resíduo quadrático. 

\item[•] Como $325=5^2\cdot13$  resíduo quadrático de $M$ então  $\mathbf{\textcolor{OliveGreen}{13}}$ é resíduo quadrático. 

\item[•] Como $13$ e $-5\cdot11\cdot13$ são resíduos quadráticos então  $-5\cdot11\cdot13^2$ também é resíduo quadrático o que implica que $\mathbf{\textcolor{OliveGreen}{-5\cdot11}}$ seja também um resíduo quadrático. 

\end{itemize}

\n Seja $C$ qualquer  classe diferente da classe principal de formas de um determinante negativo $-M$ ou mais geral $-kM$ e seja seu período $2C,3C,\cdots.$\\
As classes $2C,4C,\cdots$  pertencem ao género principal; $3C,5C,\cdots$ pertencem ao mesmo género que $C$. Se portanto $(a,b,c) $ é a forma mais simples em $C$ e $(a',b',c')$
 uma forma em alguma classe do período , digamos $nC$, ou $a'$ ou $aa'$ será um resíduo de $M$ segundo $n$ seja par ou ímpar(no primeiro caso $c'$  será também resíduo quadrático, no ultimo caso $ac',ca'\  e\  cc' $ serão resíduos quadráticos). \\
 
\n O calculo do período, isto é, das formas mais simples nas suas classes, é fácil quando $a$ é muito pequeno, especialmente quando é igual a 3, o que é sempre possível quando $kM\equiv 2\ mod\ 3$
\begin{exemplo}
Inicio do período da classe que contém a forma $(3,1,332444).$
\begin{itemize}
\item[•] $C(3,1,332444).$
\item[•] Calculemos $2C=(A,B,C)$ composição da forma $(a,b,c)=(3,1,332444)$ com $(a',b',c')=(3,1,332444)$.\\
$h=3$, $b+b'=2$, $mdc(a,a')=3$ , $(3,2)=1$\\
$A=3^2=9$, $B\equiv 997333\ mod\ 9 \Rightarrow B=-2.$\\
$C= \frac{-2^2+997331}{9}=110815$.\\
Logo, $\mathbf{2C=(9,-2,110815).}$

\item[•] Composição das formas $(a,b,c)=(9,-2,110815)$ e $(a',b',c')=(3,1,332444)$.\\
$h=3$, $b+b'=-1$, $mdc(3,3^2)=3$\\
$A=3^3=27$, $B\equiv 1994672\ mod\ 27 \Rightarrow B=7.$\\
$C= \frac{7^2+997331}{27}=36940$.\\
Logo, $\mathbf{3C=(27,7,36940)}.$
\item[•]  $\mathbf{4C=(81,34,12327)}.$
\item[•]  $\mathbf{5C=(243,34,4109)}.$
\item[•]  $\mathbf{6C=(729,-209,1428)}.$
\item[•]Composição das formas $(a,b,c)=(3,1,332444)$ e $(a',b',c')=((729,-209,1428)$.\\
$h=3$, $b+b'=-208$, $mdc(3,3^6)=3$\\
$A=3^7=2187$, $B\equiv 2081815\ mod\ 2187 \Rightarrow B=209.$\\
$C= \frac{209^2+997331}{2187}=476$.\\
Logo, $7C=(2187,209,476).$\\ 
A forma $(2187,209,476)$ é propriamente equivalente à forma $\mathbf{(476,209,2187)}.$
\item[•] $8C=(1428,685,1027).$ é propriamente equivalente à forma $\mathbf{(1027,342,1085)}.$
\item[•] $9C=(3081,1369,932).$ é propriamente equivalente à forma $\mathbf{(932,-437,1275)}.$
\item[•] $10C=(2796,2359,2347).$ é propriamente equivalente à forma $\mathbf{(425,12,2347)}.$
\end{itemize}

Assim, temos as seguintes formas \\

$\begin{array}{rrrr}
C(&3,&1,&332444)\\
2C(&9,&-2,&110815)\\
3C(&27,&7,&36940)\\
4C(&81,&34,&12327)\\
5C(&243,&34,&4109)\\
6C(&729,&-209,&1428)\\
7C(&476,&209,&2187)\\
8C(&1027,&342,&1085)\\
9C(&932,&-437,&1275)\\
10C(&425,&12,&2347)\\
\end{array}$\\

\n Das formas $C$ e $7C$ temos o resíduo $3\cdot 476 = 3 \cdot 2^2 \cdot 7 \cdot 17$. Logo, $\mathbf{\textcolor{OliveGreen}{3\cdot 7 \cdot 17}}$ também é resíduo quadrático $mod\ 997331.$\\

\n Da forma $8C$ obtemos os resíduos $1027 = \mathbf{\textcolor{OliveGreen}{13 \cdot 79}}$ e $1085= \mathbf{\textcolor{OliveGreen}{5\cdot 7 \cdot 31}} $ $mod\ 997331.$\\

\n Da forma $10C$ temos o resíduo $425= 5^2 \cdot 17$. logo daqui obtemos o resíduo quadrático $\mathbf{\textcolor{OliveGreen}{17}}$ $mod\ 997331.$
\end{exemplo}

\n Considerando os resíduos do exemplo \ref{e2.1}: $37,2\cdot 5\cdot67,-17\cdot83,-5\cdot11\cdot13,-17\cdot 53,-2\cdot 59,-2\cdot3\cdot 17, 13, -5\cdot11 .$, então 
\begin{itemize}
\item[•] De $ -2\cdot3\cdot 17$  e $17$ obtemos o resíduo $\mathbf{\textcolor{OliveGreen}{-2\cdot 3}}.$
\item[•] De $ -2\cdot3\cdot 17$  e $3\cdot 7 \cdot 17$ obtemos o resíduo  $\mathbf{\textcolor{OliveGreen}{-2\cdot 7}}.$
\item[•] De $ -17\cdot 53$  e $17$ obtemos o resíduo  $\mathbf{\textcolor{OliveGreen}{-53}}.$
\item[•] De $ 13\cdot 79$  e $13 $ obtemos o resíduo  $\mathbf{\textcolor{OliveGreen}{79}}.$
\item[•] De $ -17\cdot 83$  e $17 $ obtemos o resíduo  $\mathbf{\textcolor{OliveGreen}{-83}}.$
\item[•] De $ 5 \cdot 7 \cdot 31$  e $-2\cdot 7$ obtemos o resíduo  $\mathbf{\textcolor{OliveGreen}{-2\cdot 5 \cdot 31}}.$
\end{itemize}

\n Em resumo temos até o momento  os seguintes resíduos quadráticos: 
\[-6,13,-14,17,37,-53,-55,79,-83,-118,-310,670.\]

\n Agora, usando o método de Gauss, já estudado, obtemos outros resíduos quadráticos $mod\ 997331.$\\

\vspace*{.3cm}

\hspace*{6cm} \textbf{Resíduos quadráticos obtidos}
\begin{itemize}

\item[•] $\begin{array}{lrrrrrrrrrrrrr}
 997331 = 999^2- 2 \cdot 5 \cdot 67  & \hspace*{5cm} &  2 \cdot 5 \cdot 67.\end{array}$ 

\item[•] $\begin{array}{lrrrrrrrrrrrrr}997331 = 994^2+ 5 \cdot 11 \cdot 13^2 &\hspace*{4.8cm} & -5 \cdot 11 .\end{array}$ 

\item[•] $\begin{array}{lrrrrrrrrrrrrr}997331 = 2\cdot 706^2 + 3 \cdot 17 \cdot 3^2 &\hspace*{2.5cm} & -2\cdot 3 \cdot 17 = \mathbf{\textcolor{OliveGreen}{-102}}.\end{array}$

\item[•] $\begin{array}{rrrrrrrrrrrrr}997331 = 3\cdot575^2- 11 \cdot 31 \cdot 4^2 & \hspace*{3.8cm}&  -3 \cdot 11 \cdot 31.\end{array}$ 

\item[•] $\begin{array}{rrrrrrrrrrrrr}997331 = 3 \cdot 577^2- 7 \cdot 13 \cdot 4^2 & \hspace*{4.5cm}&3 \cdot 7 \cdot 13.\end{array}$ 

\item[•] $\begin{array}{rrrrrrrrrrrrr}997331 = 3 \cdot 578^2- 7 \cdot 19 \cdot 37& \hspace*{3.8cm}& 3 \cdot 7 \cdot 19 \cdot 37.\end{array}$ 

\item[•] $\begin{array}{rrrrrrrrrrrrr}997331 =11 \cdot 299^2- 2 \cdot 3 \cdot 5 \cdot 29 \cdot 4^2 & \hspace*{1.8cm}& -11 \cdot 2 \cdot 3 \cdot 5 \cdot 29.\end{array}$ 

\item[•] $\begin{array}{rrrrrrrrrrrrr}997331 =11 \cdot 301^2 + 5 \cdot 12^2& \hspace*{5cm}& -11 \cdot 5.\end{array}$ 

\end{itemize}

\n Dos resíduos $3 \cdot 7 \cdot 13$ e $3 \cdot 7 \cdot 19 \cdot 37$ obtemos o resíduo $13  \cdot 19 \cdot 37$. Como $13$ e $37$ já sabemos que são resíduos então obtemos o resíduo $\mathbf{\textcolor{OliveGreen}{19}}.$

\n Dos resíduos $ -11 \cdot 2 \cdot 3 \cdot 5 \cdot 29$ e $-11 \cdot 5$ obtemos o resíduo $2\cdot 3 \cdot 29$. Como -6  já sabemos que é resíduo então obtemos o resíduo $\mathbf{\textcolor{OliveGreen}{-29}}.$

\n Logo, dos métodos anteriores obtemos os resíduos

\[-6,13,-14,17,19,-29,37,-53,-55,79,-83,-102-118,-310,670.\]

\n Se algum dos números  anteriores for um não-resíduo de algum dos primos menor que $998$ então este primo  é eliminado como um possível divisor de 997331.\\
Ao realizar a verificação anterior concluímos que o único primo para o qual todos os números anteriores são resíduos quadráticos é 127. \\

\n Logo 127 é um divisor de 997331. De fato, $997331=127 \cdot 7853.$




\end{document}